\documentclass[12pt, reqno]{amsart}
\usepackage{amsmath}

\usepackage{amsfonts}
\usepackage{amssymb}
\usepackage{comment}
\usepackage{graphicx}
\usepackage{xcolor}
\usepackage{enumitem}
\usepackage[margin=1in,footskip=0.25in]{geometry}
\providecommand{\U}[1]{\protect\rule{.1in}{.1in}}
\newtheorem{theorem}{Theorem}

\newtheorem{corollary}[theorem]{Corollary}

\newtheorem{definition}[theorem]{Definition}
\newtheorem{defrem}[theorem]{Definition and Remark}

\newtheorem{lemma}[theorem]{Lemma}

\newtheorem{proposition}[theorem]{Proposition}
\newtheorem{remark}[theorem]{Remark}

\begin{document}

\title[Coorbit spaces over local fields]{Wavelet Coorbit Spaces over Local Fields}


\author{Kumar Abhinav, Hartmut F\"uhr, Qaiser Jahan}
\address{School of Mathematical and Statistical Sciences, IIT Mandi,
India}

\email{d19021@students.iitmandi.ac.in}

\address{Lehrstuhl f\"ur Geometrie und Analysis, RWTH Aachen University, D-52056 Aachen,
Germany}

\email{fuehr@mathga.rwth-aachen.de}

\address{School of Mathematical and Statistical Sciences, IIT Mandi,
India}

\email{qaiser@iitmandi.ac.in}

%
%
%
%




\begin{abstract}
This paper studies wavelet coorbit spaces on disconnected local fields $K$, associated to the quasi-regular representation of $G = K \rtimes K^*$ acting on $L^2(K)$. We show that coorbit space theory applies in this context, and identify the homogeneous Besov spaces $\dot{B}_{\alpha,s,t}(K)$ as coorbit spaces. We identify a particularly convenient space $\mathcal{S}_0(K)$ of wavelets that give rise to tight wavelet frames via the action of suitable, easily determined discrete subsets $R \subset G$, and show that the resulting wavelet expansions converge simultaneously in the whole range of coorbit spaces. For orthonormal wavelet bases constructed from elements of $\mathcal{S}_0(K)$, the associated wavelet bases turn out to be unconditional bases for all coorbit spaces. We give explicit constructions of tight wavelet frames and wavelet orthonormal bases to which our results apply. 
\end{abstract}




\maketitle

\section{Introduction}\label{sec1}

This paper is primarily concerned with the construction of discrete wavelet systems in local fields and their use in the characterization and expansion of a wide variety of function spaces. The construction of wavelet systems on local fields has been an active field since last three decades. J. J. Benedetto and R. L. Benedetto had developed a wavelet theory for general local fields in \cite{bib9, bib10}. 

Essentially, the local fields are of two types (excluding the connected local fields which are $\mathbb{R}$ and $\mathbb{C}$). The characteristic zero local field includes the $p-$adic field $\mathbb{Q}_p$ and of positive characteristic include the $p-$series field, Cantor dyadic group, and Vilenkin groups. The Characteristic of the field plays a crucial role in the construction of wavelets. Their wavelet theory are quite different even though their structure and metrics are similar for both the local fields.

In 2002, Kozyrev \cite{MR1918846} constructed the first wavelet basis for $L^2(\mathbb{Q}_p)$. The general scheme for the construction of wavelets in $\mathbb{Q}_p$ using multiresolution analysis (MRA) was given by Shelkovich and Skopina in \cite{MR2511868}. Khrennikov et al. \cite{MR2558153} constructed a number of scaling functions, but later on in \cite{MR2673705} it was shown that they are all associated with the same Haar MRA. For the Cantor dyadic group, we refer to \cite{MR1373159}. Farkov \cite{Far} constructed wavelets in Vilenkin groups. Behera and Jahan have developed a theory of wavelet on local field of positive characteristic \cite{mabook}.

Concurrently, the theory of homogeneous and inhomogeneous Besov spaces on these groups was developed, essentially by analogy to the Euclidean case (see e.g. \cite{bib8}).

A fundamental connection between wavelets and homogeneous Besov spaces in the euclidean setting is provided by coorbit space theory, which develops tools for the definition and study of function spaces via the decay of expansion coefficients. The parallels to the Euclidean case are clear enough to allow a straightforward adaptation of the arguments. Just as in the Euclidean case, the wavelet systems arise as (subsets of) orbits of the so-called quasi-regular representation $\pi$ of affine group on local fields, $K$, i.e., the semidirect product group $G = \mathbb{K} \rtimes K^*$ acting on $L^2(K)$, and the applicability of coorbit space theory hinges on verifying certain integrability properties of this representation.  To our knowledge, this connection has not yet been fully explored.

This paper serves to provide the missing link, and to illustrate that the continuous wavelet transform and its associated coorbit spaces are indeed very useful both for the construction of tight wavelet frames on local fields, and for the characterization of a large variety of function spaces. A second, somewhat related reason is the topogical structure of $K$ and its affine group $G = \mathbb{K} \rtimes K^*$. Both have a neighborhood basis at unity consisting of open, compact subgroups, which allows to reduce sampling theorems to simple cases of integration over disjoint cosets. 

As will emerge below, wavelet coorbit theory over local fields is in many respects cleaner and much more powerful than its euclidean counterpart. The main reason for this is the availability of an easily accessible reservoir of fundamental building blocks (wavelets), for which all relevant properties needed for the proper functioning of coorbit theory are actually easy to establish. The building blocks in question are the usual test functions on $\mathcal{S}(K)$, namely compactly supported functions that are locally constant. The subspace $\mathcal{S}_0(K)$ of test functions with vanishing integral turns out to have many advantageous properties that lead to much stronger statements than coorbit theory readily supplies in the euclidean case.

A quick summary of the main contributions of this paper is as follows:
\begin{itemize}
\item Given any $0 \not= g \in \mathcal{S}_0(K)$, there exists a subset $R \subset K \rtimes K^*$, explicitly computable from $g$, such that the wavelet system $(\pi(r) g)_{r \in R} \subset L^2(K)$ is a tight frame (Theorem \ref{thm:tight_frame}).
\item Given a tight wavelet frame associated to a wavelet $g \in \mathcal{S}_0(K)$, the functions $f \in S_0(K)$ are characterized by the fact that only finitely many wavelet coefficients are nonzero (Corollary \ref{cor:wavelet_char_S0}).
\item Every $\varphi$ in the dual space $\mathcal{S}_0'(K)$ has a wavelet expansion converging in the weak-$*$-topology (Proposition \ref{prop:expansion_S_strich}).
\item The quasi-regular representation $\pi$ is integrable, allowing the definition of a large variety of coorbit spaces (Corollary \ref{cor:coorbit_appl}).
\item The homogeneous Besov spaces $\dot{B}_{s,t,\alpha}(K)$, for $1 \le s,t \le \infty$, are coorbit spaces. (Theorem \ref{thm:besov_coorbit}).
\item A tight wavelet frame of $L^2(K)$, generated from a finite set of mother wavelets in $\mathcal{S}_0(K)$, is a Banach frame, simultaneously for all wavelet coorbit spaces. Given an element $f \in Co(Y)$, the wavelet expansion of $f$ converges in the norm of $Co(Y)$ (Theorem \ref{thm:disc_char_coorbit}).
\item An orthonormal basis, generated from a finite set of mother wavelets in $\mathcal{S}_0(K)$, is an unconditional basis in all wavelet coorbit spaces (Theorem \ref{thm:disc_char_coorbit}).
\end{itemize}

\section{Preliminaries on Local Fields}\label{sec2}

In this section, we list some basic terminology and properties of local fields. For proofs and further details on this topic, we refer \cite{mabook,taible}.

Let $\mathit{K}$ be an algebraic field which is also a topological space. If additive group $\mathit{K}^+$ and multiplicative group $\mathit{K}^*$ of $\mathit{K}$ are locally compact abelian group then $\mathit{K}$ is called a locally compact field or a \emph{local field}.\par
Any field $\mathit{K}$ endowed with discrete topology is a locally compact field. Also if  $\mathit{K}$ is non-discrete and connected then  $\mathit{K}$ is either $\mathbb{R}$ or $\mathbb{C}$. If $\mathit{K}$ is not connected then $\mathit{K}$ is totally disconnected. So, we mean by a local field $\mathit{K}$, a locally compact, non-discrete, totally disconnected field. \par
Let $\mathit{K}$ be a local field. Since $\mathit{K}^+$ is a locally compact abelian group, we choose a Haar measure $dx$ for $\mathit{K}^+$. For $\alpha \neq 0$ and $\alpha \in \mathit{K}$, $d(\alpha x)$ is also a Haar measure. Let $d(\alpha x) = |\alpha|dx$, we call $|\alpha|$ the absolute value or \emph{valuation} of $\alpha$. Also let $|0|=0$. The absolute value has the following properties:- \par
\begin{itemize}
    \item $|x|=0$ if and only if $x=0$.
    \item $|xy|=|x||y|$ for all $x,y \in \mathit{K}$.
    \item $|x+y| \leq \,\max\,\{|x|,|y|\}$ for all $x,y \in \mathit{K}$.
\end{itemize}
The last property is called \emph{ultrametric inequality}. It also follows that if $|x| \neq |y|$ then $|x+y|=\,max\,\{|x|,|y|\}$. Haar measure on $\mathit{K}^*$ is given by $\dfrac{dx}{|x|}$ where $dx$ is a Haar measure on $\mathit{K}^+$.\par
We define $\mathfrak{D}=\{x \in \mathit{K}\,\,:\,\,|x|\leq 1\}$ and $\mathfrak{D^*}=\{x \in \mathit{K}\,\,:\,\,|x| = 1\}$. $\mathfrak{D}$ is called the \emph{ring of integers}. It is the unique maximal compact subring of $\mathit{K}$. Also define $\mathfrak{P}=\{x \in \mathit{K}\,\,:\,\,|x| < 1\}$. $\mathfrak{P}$ is called the prime ideal in $\mathit{K}$. It is the unique maximal ideal in $\mathfrak{D}$. \par

Since $\mathit{K}$ is totally disconnected, the set of absolute values is of the form $\{s^k \,:\,k \in \mathbb{Z}\} \cup \{0\}$ for some $s > 0$. In fact it is known that $s = p^\ell$, for a suitable prime number $p$ and $\ell \in \mathbb{N}$. By compactness there is an element in $\mathfrak{P}$ of maximum absolute value. Let $\mathfrak{p}$ be such an element in $\mathfrak{P}$. Then $\mathfrak{p}$ is called a prime element in $\mathit{K}$ and $\mathfrak{P} = (\mathfrak{p}) = \mathfrak{p}\mathfrak{D}$ as an ideal in $\mathfrak{D}$.\par
$\mathfrak{D}$ is compact and open and so $\mathfrak{P}$ is also compact and open. Therefore $\mathfrak{D}/\mathfrak{P}$ is isomorphic to a finite field $GF(q)$ where $q=p^c$ for some prime number $p$ and $c \in \mathbb{N}$. We refer \cite{mabook} for the proof. \par
For a measurable set $S \subset \mathit{K}$, let $|S| = \int_{\mathit{K}} \mathcal{X}_S\, dx$, where $\mathcal{X}_S$ is the characteristic function of $S$, $dx$ is Haar measure normalized so that $|\mathfrak{D}|=1$. Then $|\mathfrak{P}|=q^{-1}$, and $|\mathfrak{p}|=q^{-1}$. It follows that for any $x \neq 0$, $x \in \mathit{K}$, $|x|=q^k$ for some $k \in \mathbb{Z}$.\par
We also define $\mathfrak{P}^k:= \mathfrak{p}^k\mathfrak{D} = \{x \in \mathit{K} \,:\,|x| \leq q^{-k}\}$  and   $ \mathfrak{D}_k^*:= \mathfrak{p}^k\mathfrak{D} = \{x \in \mathit{K} \,:\,|x| = q^{-k}\} , k \in \mathbb{Z}$. $\mathfrak{P}^K$ is called a fractional ideal. $\mathfrak{P}^k$ is compact and open and is a subgroup of $\mathit{K}^+$ for all $k \in \mathbb{Z}$.\par

Let $\mathit{K}$ be a local field then there is a non-trivial, unitary, continuous character $\chi$ on $\mathit{K}^+$. The existence of such a character is a consequence of Pontryagin duality theorem. If $\chi$ is a nontrivial character on $\mathit{K}^+$ then there is a $k \in \mathbb{Z}$ such that $\chi$ is trivial on $\mathfrak{P}^k$. Then $\chi$ is also constant on the cosets of $\mathfrak{P}^k$ in $\mathit{K}^+$. Observe that $\chi$ do identify $\mathit{K}$ with its dual, i.e., $\mathit{K}^+ \cong \hat{\mathit{K}}^+$ with the correspondence $\lambda \longleftrightarrow \chi_{\lambda}$ where $\chi_{\lambda}(x) := \chi (\lambda x)$. For the proof, see \cite{taible}. Throughout the paper, we will employ a character $\chi$ that has $\mathfrak{D}$ as its kernel to identify $K$ with $\hat{K}$.\\

\begin{definition}
    If $f \in L^1(\mathit{K})$, the Fourier transform of $f$ i.e., $\hat{f}$ is defined as
    $$\hat{f}(x) = \int_{\mathit{K}} f(\xi) \overline{\chi_x(\xi)}\, d\xi = \int_{\mathit{K}} f(\xi) \chi(-x \xi)\, d\xi $$
\end{definition}
Fourier transform of $f$ has the following properties:-
\begin{itemize}
    \item The map $f \mapsto \hat{f}$ is a bounded linear transformation of $L^1$ into $L^{\infty}$ i.e., $\|\hat{f}\|_{\infty}\leq \|f\|_1$.
    \item  If $f \in L^1$ then $\hat{f}$ is uniformly continuous. 
\end{itemize}
For $k \in \mathbb{Z}$, let $\Phi_k$ be the characteristic function of $\mathfrak{P}^k$. Then note that the characteristic function of $h + \mathfrak{P}^k$ is $\tau_h \Phi_k$ and that  $\tau_h \Phi_k$ is constant on the cosets of $\mathfrak{P}^k$.\par
Let $\mathcal{S}(\mathit{K})$ be the space of all finite linear combinations of the form $\tau_h \Phi_k, \, h \in \mathit{K}, \, k \in \mathbb{Z}$. Then $\mathcal{S}(\mathit{K})$
is an algebra of continuous functions with compact support that separate points. Consequently, $\mathcal{S}(\mathit{K})$ is dense in $C_0$ as well as in $L^p,\, 1 \leq p < \infty $.

Alternatively, we can describe $\mathcal{S}(\mathit{K})$ as follows:
$\phi \in \mathcal{S}(\mathit{K})$ if and only if there are integers $k,l$ such that $\phi$ is constant on cosets of $\mathfrak{P}^k$ and is supported on $\mathfrak{P}^l$. 
We define this subspace as
$$\mathcal{S}^k_l:= \{\phi \in \mathcal{S}(K)\,:\, \phi \, \text{is supported on} \,\mathfrak{P}^l\, \text{ and is constant on the cosets of}\,\mathfrak{P}^k\}.$$
Clearly this space is nontrivial only if $l \le k$.
By the above remarks we have $\mathcal{S}(K) = \bigcup_{l,k \in \mathbb{Z}} \mathcal{S}^k_l$.

If $\phi \in \mathcal{S}_{l}^k$, then $\hat{\phi} \in \mathcal{S}_{-k}^{-l}$, by \cite[Theorem 3.2]{taible}. In particular, $\mathcal{S}(K)$ is invariant under the Fourier transform. In fact, the Fourier transform is a homeomorphism of $\mathcal{S}(\mathit{K})$ onto $\mathcal{S}(\mathit{K})$, by Theorem I.(3.2) of \cite{taible}.

 $\mathcal{S}(K$) is topologized as follows: For each $f \in \mathcal{S}(K)$, a neighborhood base is given  by the sets $U \subset \mathcal{S}(K)$,
\[
U = \{ g \in \mathcal{S}_l^k : \| f - g \|_\infty < \epsilon \} 
\] with the additional assumptions $f \in \mathcal{S}_l^k$ and $\epsilon>0$. This makes $\mathcal{S}(K)$ a topological vector space. It is easy to see that the containment $\mathcal{S}(\mathit{K}) \subset L^p(G)$ is continuous.

We call $\mathcal{S}(\mathit{K})$, the space of test functions. $\mathcal{S}(\mathit{K})$ is a complete, separable locally convex space. A subspace of particular relevance is given by
\[                                                                                                                                                \mathcal{S}_0 = \mathcal{S}_0(\mathit{K}) = \left\{ f \in \mathcal{S}(\mathit{K}) : \int_K f(x) dx = 0 \right\}~.                                                                                                                                                                                        \] The continuity of the embedding $\mathcal{S}(\mathit{K}) \subset L^1(G)$ ensures that $\mathcal{S}_0(K) \subset \mathcal{S}(\mathit{K})$ is a closed subspace. We denote its dual by $\mathcal{S}_0'(K)$.

\section{Assumptions for Coorbit Theory}\label{sec3}

This section summarizes the basic assumptions, terminology and certain technical details required to define coorbit spaces. \par
Let $G$ denotes a $\sigma$-compact, locally compact group. Let $dx$ denotes the left Haar measure of $G$ and $\Delta$ the modular function. We define the following operations on the functions on $G$:-
\begin{itemize}
    \item Left translation by $x \in G$ :   $L_xf(y) = f(x^{-1}y)$
    \item Right translation by $x \in G$ :  $R_xf(y) = f(yx)$
    \item Involution : $f^{\nabla}(y) = \overline{f(y^{-1})}$\\
\end{itemize}

\begin{definition}
 Given a locally compact group $G$, a (submultiplicative)\textbf{weight} on $G$ is a Borel measurable function $w: G \to \mathbb{R}^+$ that is bounded on compact sets and satisfies
 \[
\forall x,y \in G~:~ w(xy) \le w(x) w(y)~.
 \]

 Given two weights $w$ and $v$, we write $w \preceq v$ if $w(g) \le C v(g)$ holds for all $g \in G$ and some constant $C>0$. We let $w \sim v$ if $w \preceq v$ and $v \preceq w$.

 We call a weight $w$ on $G = \mathit{K} \rtimes \mathit{K}^*$ \textbf{separable} if $w(x,h) = w(x,1) w(0,h)$ holds.
\end{definition}

\begin{definition}
    We call $Y$ a \textbf{solid Banach function space} if
   \begin{enumerate}
       \item $Y$ is continuously embedded into locally integrable functions $L^1_{Loc}(G)$ on $G$.
       \item $Y$ is solid i.e., if $|f| \leq |g|$ a.e. and $g \in Y$ then $f \in Y$ and $\left\|f\right\|_Y \leq \left\|g\right\|_Y$.\\
   \end{enumerate}
\end{definition}

\begin{definition}
    Let $G$ be a locally compact Hausdorff group and choose $Q \subset G$  measurable, precompact. Define the \textbf{ maximal function} of $f : G \to \mathbb{C}$  as
    $$M_Qf : G \to [0, \infty]\,\, ,\,\,x \mapsto {\rm ess sup}_{y \in xQ} |f(y)|.$$
    For a solid Banach function space $(Y, \left\|.\right\|_Y)$, we define the \textbf {Weiner amalgam space} as 
    $$ W_Q(L^{\infty},Y) :=\{ f \in L^{\infty}_{Loc} (G)\,:\, M_Qf \in Y \}$$
    with $$\left\|f\right\|_{W_Q(Y)} := \left\|f\right\|_{W_Q(L^{\infty},Y)}:= \left\|M_Qf\right\|_Y $$
    for $f\in W_Q(L^{\infty},Y)$.
    
    For non-abelian groups, \textbf{right-sided maximal function} of $f$ is defined as
    $$M_Q^Rf : G \to [0, \infty]\,\,,\,\,x \mapsto {\rm ess sup}_{y \in Qx} |f(y)|.$$
    Then we define \textbf{right-sided Wiener amalgam space} as
     $$W_Q^R(L^{\infty},Y) :=\{ f \in L^{\infty}_{Loc} (G)\,:\, M^R_Qf \in Y \}$$
    with $$\left\|f\right\|_{W_Q^R(L^{\infty},Y)}:= \left\|M_Q^Rf\right\|_Y .$$
\end{definition}
 Both $W_Q(L^{\infty},Y)$ and $W_Q^R(L^{\infty},Y)$ embed continuously into $Y$ and are independent of the choice of $Q$.\\

\begin{definition}
    Let $\pi$ be an irreducible, unitary, continuous representation $\pi$ of $G$ on a Hilbert space $\mathcal{H}$.  For $f,g \in \mathcal{H}$, we define the \textbf{(generalized) Wavelet Transform} of $f$ (with respect to g) as
    $$W_gf:G \to \mathbb{C}\,\,,\,\, x \mapsto \langle f|\pi(x)g\rangle_{\mathcal{H}}.$$
    It follows from Cauchy-Schwarz inequality and strong continuity of $\pi$ that $W_g$ is always a bounded, continuous function. The representation $\pi$ is called \textbf{square-integrable} if there exists some admissible $g \in \mathcal{H} \setminus \{0\}$. Here $g \in \mathcal{H} $ is called admissible if $W_gg \in L^2(G)$. The representation $\pi$ is called \textbf{integrable} if $W_gg \in L^1(G)$ for some $g \in \mathcal{H} \setminus \{0\}$.
\end{definition}
 
 We next give a list of the conditions and definitions required for the formulation of coorbit theory (cf.\cite{bib2},\cite{bib3}). In the following, we assume that $G$ is a $\sigma$-compact, locally compact Hausdorff group. \par

 \begin{enumerate}[label=\alph*)]
     \item $\pi : G \to \mathcal{U}(\mathcal{H)}$ is a strongly continuous, unitary, irreducible, integrable (and hence square integrable) representation of $G$ for some non-trivial Hilbert space $\mathcal{H}$.
     \item Let $(Y, \left\|\cdot\right\|_Y)$ be a solid- Banach function space which is invariant under left/right translations.
     \item Let $v: G \to (0,\infty)$ be a fixed measurable, submultiplicative weight with $$v(x) \geq \, \max\, \{\left\|L_{x^{\pm 1} }\right\|_{Y \to Y},\,\left\|R_x\right\|_{Y \to Y},\,\left\|R_{x^{-1}} \right\|_{Y \to Y}\Delta(x^{-1})\}$$
     and also $$v(x) = v(x^{-1}) \Delta(x^{-1}).$$
     In particular, $v(x) \geq 1,\,\,\, \left\|f\right\|_{L_v^1} =\left\|f^{\nabla}\right\|_{L_v^1} $ and $Y \ast L_v^1 \subseteq Y$. Such a weight is called \textbf{control weight} for $Y$.
     \item We define the set of \textbf{analyzing vectors}
     $$\mathcal{G}_v = \{ g \in \mathcal{H}\,\,: W_g g \in L_v^1\}$$
     and we assume that $\mathcal{G}_v$ is non trivial.
     \item For a fixed non trivial vector $g \in \mathcal{G}_v$, we define the space of test vectors     $$\mathcal{H}_v^1 :=\,\{f \in \mathcal{H}\,:\,W_g f \in L_v^1\}$$
     with norm
     $$\left\|f\right\|_{\mathcal{H}_v^1}:=\,\left\|f\right\|_{\mathcal{H}_v^1}:=\,\left\|W_gf\right\|_{L_v^1}.$$
 \end{enumerate} 
 Let $\mathcal{H}_v^{1,\sim}$ be the space of all continuous, conjugate-linear functionals on $\mathcal{H}_v^1$. The inner product on $\mathcal{H} \times \mathcal{H}$ extends to a sesquilinear form on $\mathcal{H}_v^1 \times \mathcal{H}_v^{1,\sim} $. So for $g \in \mathcal{H}_v^1$ and $f \in \mathcal{H}_v^{1,\sim}$, $W_g f (x) = \langle f \,|\,\pi (x)g\rangle$ is well defined. \par
 Now, for a fixed analyzing vector $g \in \mathcal{G}_v$, we define the coorbit space with respect to $Y$ and $\pi$
 $$Co(Y):=\,Co_g(Y):=\, \{f\in \mathcal{H}_v^{1,\sim}\,:\, W_g f \in Y\}$$
 with the norm 
 $$\left\|f\right\|_{Co(Y)}:=\,\left\|f\right\|_{Co_g(Y)}:=\,\left\|W_gf\right\|_{Y}.$$
 One of the most important properties of coorbit spaces is that $Co_g(Y)$ is independent of the choice of analyzing vector $g \in \mathcal{G}_v \setminus \{0\}$ \cite[Lemma 4.2]{bib2}. We define another class of vectors called \textbf{better vectors} as
 $$\mathcal{B}_v:=\,\{g\in \mathcal{H}\,:\,W_g g \in W^R(L^{\infty},Y)\}.$$

Better vectors are sometimes called ``atoms" as they are used to obtain atomic decomposition of coorbit spaces. They are also used to define coorbit spaces in more generality (as for quasi-Banach function spaces). Since $W_Q^R(L^{\infty},Y)$ has continuous embedding into $Y$, $\mathcal{B}_v \subset \mathcal{G}_v$.

 \section{Applicability of Coorbit Theory for the Spaces $L_v^{p,q}(\mathit{K})$}\label{sec4}

In this section, we lay the foundation for the study of wavelet coorbit spaces on a local field. We assume $G = \mathit{K} \rtimes \mathit{K}^*$ where $\mathit{K}$ is a local field and $\mathit{K}^* \subset \mathit{K}$ is the multiplicative group of $\mathit{K}$. Multiplication on $G$ is defined as
$$(x,h)\cdot(y,k) = (hy + x, hk)$$ for all $(x,h), (y,k) \in G$. So $G$ is a locally compact Hausdorff group with Haar integral given by
$$\int_G f(x,h)\,d(x,h) = \int_{K^*}\int_K f(x,h) \, dx \frac{dh}{|h|}$$
where $dh$ denotes the Haar measure on $\mathit{K}^*$. The modular function of $G$ is 
$$\Delta_G (x,h) = \frac{\Delta_{\mathit{K}^*}}{|h|}= |h|^{-1}$$ as $\mathit{K}^*$ is abelian and hence unimodular. The quasi-regular representation of $G$ is given as $$\pi : G \to \mathcal{U}(L^2(\mathit{K}))\,\,,\,\,(x,h) \mapsto L_x \Delta_h$$
where $\Delta_hf(y) := |h|^{-1/2}f(h^{-1}y)$, $\pi$ is a unitary and strongly continuous representation of $G$. It is easy to adapt the argument for the real case to show that $\pi$ is \textit{irreducible}, and in Lemma \ref{lem:g_adm} below it is shown that $\pi$ is \textit{square-integrable} in our case.\\
$\mathit{K}^*$ acts on $\mathit{K}$ by a topological automorphism i.e.,
$$h(a) = h\cdot a.$$
The dual action $\mathit{K}^*$ on $\hat{\mathit{K}} \cong \mathit{K}$ is given by
$$p_{\xi} : \mathit{K}^* \times \hat{\mathit{K}} \to \hat{\mathit{K}}\,\,,\,\,(h,\xi) \mapsto h^{-1}\cdot \xi.$$
 We define the \textbf{weighted mixed Lebesgue space} $L_w^{p,q}(\mathit{K})$ as
$$L_w^{p,q}(\mathit{K}):= \{f:G \to \mathbb{C} \,:\, f \,\text{ measurable and}\, \left\|f\right\|_{L_w^{p,q}(\mathit{K})} <  \infty\}$$
with $$\left\|f\right\|_{L_w^{p,q}(\mathit{K})}:=\left(\int_{\mathit{K}^*} \|w(\cdot,h)\cdot f(\cdot,h)\|_{L^p(\mathit{K})}^q\frac{dh}{|h|} \right)^{1/q} $$
for $q \in [1,\infty)$ and with
$$\left\|f\right\|_{L_w^{p,q}(\mathit{K})}:= {\rm ess~sup}_{h \in \mathit{K}} \|w(\cdot,h)\cdot f(\cdot,h)\|_{L^p(\mathit{K})} $$ for $q = \infty$.

We will only consider weights that are continuous, submultiplicative i.e., a weight $v$ such that $v(x_1x_2) \leq v(x_1)v(x_2)$, for any $x_1, x_2 \in G$. Furthermore, all weights will be assumed to be of moderate growth i.e.,
$$v(x,h) \leq (1+|x|)^sw(h)$$
for some suitable weight $w$ and $s \geq 0$.

The following lemma provide explicit estimates for the translation norms of $L_v^{p,q}(\mathit{K})$, since these spaces will be the coefficient spaces associated to the Besov spaces. For a proof, we refer to \cite[Theorem 2.2.19]{bib5}.

\begin{lemma}{\label{lroprnorm}}
    Let $G$ be a locally compact group and assume that $(Y, \left\|.\right\|_Y)$ is a left/right invariant solid function space on $G$. Let $v:G \to (0,\infty)$ be a weight, then $Y_v$ is also left/right invariant with 
    $$\left\|L_x\right\|_{Y_v \to Y_v} \leq \left\|L_x\right\|_{Y \to Y}.v(x) $$
    $$\left\|R_x\right\|_{Y_v \to Y_v} \leq \left\|R_x\right\|_{Y \to Y}.v(x^{-1}) $$
    for all $x\in G$. 
\end{lemma}

We next provide explicit expressions for the control weights associated to these coefficient spaces. The argument is close to the euclidean case treated in \cite{bib5}, and is included for self-containment. 
\begin{lemma}
    Let $v:G \to (0, \infty)$ denotes a weight such that $v(x,h) \leq (1 + |x|)^s w(h)$, for a suitable weight $w$ and $s\geq 0$. Then for the space $Y = L^{p,q}_v$, there exists a control weight such that 
    $$v_0(x,h) \leq (1 + |x|)^s w_0(h)$$
for some weight $w_0(h)$ and $$max( \left\|L_{(x,h)^{\pm1}}\right\|_{Y \to Y},\left\|R_{(x,h)}\right\|_{Y \to Y},\left\|R_{(x,h)^{\pm1}}\right\|_{Y \to Y}\Delta_G(x,h)^{-1}) \leq v_0(x,h)$$
where $L_{(x,h)}, R_{(x,h)} : Y \to Y$ are the left and right translation operators.
\end{lemma}
\begin{proof}
    We start by estimating the norms of translation operators on $Y$. For $(y,k) \in G$, we have
    $$[L_{(x,h)^{-1}}f](y,k) = f( (x,h)(y,k) )=f(x +hy, hk)$$
    $$[R_{(x,h)}f](y,k) = f( (y,k)(x,h) )=f(kx+y, kh).$$
    Let us first assume that $q \in [1, \infty)$ and $F(k) := \dfrac{\left\|f(\cdot,k)\right\|^q_{L^p} }{|k|}$, then using the change of variable formula
    \begin{equation*}
        \begin{split}
            \frac{\left\|[L_{(x,h)^{-1}}f](\cdot,k)\right\|_{L_p}^q}{|k|} &= \frac{\left\|f(h \cdot +x,hk)\right\|_{L_p}^q}{|k|} = \frac{|h|^{1-q/p}\left\|f(\cdot ,hk)\right\|_{L_p}^q}{|hk|}\\
            &=|h|^{1-q/p}  L_{h^{-1}}F(k).
        \end{split}
    \end{equation*}
    Now
      \begin{equation*}
        \begin{split}
            \left\|L_{(x,h)^{-1}}f\right\|_{L_{p,q}}^q &= |h|^{1-q/p} \int_{\mathit{K}^*} L_{h^{-1}}F\, dk = |h|^{1-q/p} \int_{\mathit{K}^*} F\, dk = |h|^{1-q/p} \int_{\mathit{K}^*} \frac{\left\|f(\cdot,k)\right\|_{L_p}^q}{|k|} \, dk \\
            &=|h|^{1-q/p} \left\|f\right\|_{L^{p,q}}.
        \end{split}
    \end{equation*}
Similarly, using the translation-invariance of $\|\cdot\|_{L_p(\mathit{K})}$
$$ \frac{\left\|[R_{(x,h)}f](\cdot,k)\right\|_{L_p}^q}{|k|} = \frac{|h|\left\|f(\cdot ,kh)\right\|_{L_p}^q}{|kh|} =|h|R_{h}F(k),$$
so
$$\left\|R_{(x,h)}f\right\|_{l^{p,q}}^q = |h| \int_{\mathit{K}^*}(R_h F)(k) \,dk = |h|\left\|f\right\|^q_{L^{p,q}}.$$

For $q = \infty$, we again calculate
\begin{equation*}
        \begin{split}
            \left\|L_{(x,h)^{-1}}f\right\|_{L_{p,\infty}} &= \|h \mapsto \| \left(L_{(x,h)^{-1}}f\right)(\cdot, k)\|_{L^p}\|_{L^{\infty}(\mathit{K})}\\
            &= \|h \mapsto \| f(h\cdot + x, hk)\|_{L^p}\|_{L^{\infty}(\mathit{K})}\\
            &= |h|^{-1/p}\|h \mapsto \| f(\cdot, hk)\|_{L^p}\|_{L^{\infty}(\mathit{K})}\\
             &= |h|^{-1/p}\|g \mapsto \| f(\cdot, g)\|_{L^p}\|_{L^{\infty}(\mathit{K})}\\
             &= |h|^{-1/p} \|f\|_{L^{p,\infty}} < \infty.
        \end{split}
    \end{equation*}

    Similarly
    \begin{equation*}
        \begin{split}
            \left\|R_{(x,h)}f\right\|_{L_{p,\infty}} &= \|h \mapsto \| \left(L_{(x,h)}f\right)(\cdot, k)\|_{L^p}\|_{L^{\infty}(\mathit{K})}\\
            &= \|h \mapsto \| f(\cdot, hk)\|_{L^p}\|_{L^{\infty}(\mathit{K})}\\
            &= \|h \mapsto \| f(\cdot, h)\|_{L^p}\|_{L^{\infty}(\mathit{K})}\\
             &= \|f\|_{L^{p,\infty}}.        
        \end{split}
    \end{equation*}

Now using Lemma \ref{lroprnorm}, we get 
$$\left\|L_{(x,h)}\right\|_{Y \to Y} \leq |h|^{1/p -1/q} . v(x,h)$$
$$\left\|L_{(x,h)^{-1}}\right\|_{Y \to Y} \leq |h|^{1/q -1/p} . v((x,h)^{-1})$$

$$\left\|R_{(x,h)}\right\|_{Y \to Y} \leq |h|^{1/q} . v(x,h) = \Delta_G(x,h)^{-1/q}.v((x,h)^{-1})$$

$$\left\|R_{(x,h)^{-1}}\right\|_{Y \to Y}.\Delta_G(x,h)^{-1} \leq |h|^{1 -1/q} . v((x,h)^{-1}). $$
Now, if we choose a symmetric submultiplicative weight $v_1$ of moderate growth on $G$ such that $$v(x,h) \leq v_1(x,h)w(h)$$
then $$v_2(x,h)=v_1(x,h)\,\text{max}\,\big(\Delta_G(0,h)^{-1/q}, \Delta_G(0,h)^{1/q-1} \big)\big(w(h) + w(h^{-1})\big)\big(|h|^{1/p-1/q} + |h|^{1/q-1/p}\big)$$
is a control weight for $Y$.
\end{proof}

\begin{remark}
  In \cite{bib6}, explicit weights were constructed for the case of $Y = L^{p,q}_v (\mathbb{R}^d)$, which can be adapted to the setting $Y = L^{p,q}_v (\mathit{K})$. For example, if $v(x,h) = (1 + |x| + |h|)^s w(h)$ then $v_1(x,h)$ in the above can be taken as $v_1(x,h) = (1+|x|+|h^{-1}x|+|h^{-1}| + |h|)^s$.

\end{remark}


\section{Matrix coefficients of functions in $\mathcal{S}_l^k$}

We first clarify some conditions concerning admissibility of the wavelet. Recall that $g \in L^2(K)$ is called admissible iff 
\[
 C_g = \|W_g g \|_2^2 
\] is finite, and that in this case one has the isometry condition
\begin{equation} \label{eqn:cwt_isom}
 \| W_g f \|_2^2 = C_g \| f \|_2^2~.
\end{equation}
In the concrete setting of the quasi-regular representation, the Fourier transform allows to compute
\begin{equation} \label{eqn:adm_cond}
 \| W_g g \|_2^2 = \int_K \frac{|\widehat{g}(\xi)|^2}{|\xi|} d\xi~,
\end{equation} see also \cite{MR2130226, bib4}.

\begin{lemma} \label{lem:g_adm}
 Assume that $g \in \mathcal{S}(K)$. Then the following are equivalent:
 \begin{enumerate}
  \item[(a)] $g$ is admissible. 
  \item[(b)] $W_g g \in C_c(G)$. 
  \item[(c)] $g \in S_0(K)$.
 \end{enumerate}

\end{lemma}

\begin{proof}
 Proof of $(a) \Rightarrow (c)$:
 By assumption, $g  \in \mathcal{S}_l^k(K)$ for some $l \le k$. In particular, $g$ is integrable and therefore $\widehat{g}$ is continuous. Assuming that
 \[
  0 \not= \int_G g(x) dx = \widehat{g}(0)~, 
 \] we obtain a neighborhood $U$ of zero in $K$ and of $\epsilon>0$ such that $|\widehat{g}(\xi)|^2 \ge \epsilon>0$. On the other hand, our admissibility assumption entails
 \[
  \| W_g g \|_2^2 = \int_G |\widehat{g}(\xi)|^2 \frac{d\xi}{|\xi|} < \infty~,
 \]
and thus 
 \[
  \int_U \frac{1}{|\xi|} d\xi < \infty~.
 \] 
W.l.og. $U = \mathfrak{P}^m$ for suitable $m$, which then allows to write 
\begin{eqnarray*}
 \int_U \frac{1}{|\xi|} d\xi & = & \sum_{l=m}^\infty \int_{\mathfrak{P}^l \setminus \mathfrak{P}^{l+1}} \frac{1}{|\xi|} d\xi  \\
 & = & \sum_{l=m}^\infty \int_{{\mathfrak{P}^l}\setminus{\mathfrak{P}^{l+1}}} q^l d\xi \\
 & = & \sum_{l=m}^\infty (q-1)\cdot q^{-l-1} \cdot q^l = \infty~,
\end{eqnarray*} yielding the desired contradiction. 

For the proof of $(c) \Rightarrow (b)$ we recall that $g\in \mathcal{S}^k_l(K)$ entails $\widehat{g} \in \mathcal{S}^{-l}_{-k}(K)$. Assuming $(c)$, this implies $\widehat{g}(\xi)= 0$ for all $\xi \in \mathfrak{P}^{-l}$. In addition, we have $\widehat{g} (\xi) = 0$ whenever $|\xi| \not\in \mathfrak{P}^{-k}$, finally implying 
\[
 {\rm supp}(\widehat{g}) \subset \mathfrak{P}^{-k} \setminus \mathfrak{P}^{-l}~.
\]
As a consequence, 
\begin{equation} \label{eqn:prod_four}
 \overline{\widehat{g}(h\xi)} \cdot \widehat{g}(\xi) =0 
\end{equation} holds, for all $\xi \in K$, whenever
\[
 |h| < q^{k-l} \mbox{ or } |h| > q^{l-k}~.
\] We have $W_g g (\cdot,h) = g \ast \pi(0,h) g^*$, with $g^*(x) = \overline{g(-x)}$. Hence the convolution theorem entails via (\ref{eqn:prod_four}) that 
\[
 W_g g (\cdot,h) = 0~,\mbox{ for } |h| < q^{l-k} \mbox{ or } |h| > q^{k-l}.
\]
Furthermore, for $q^{l-k} \le |h| \le q^{k-l}$, 
\[
 {\rm supp}(W_g g) (\cdot,h) \subset \mathfrak{P}^{2l-k}~,
\] and therefore finally 
\[
 {\rm supp}(W_g g)  \subset  \mathfrak{P}^{2l-k} \times \left( \mathfrak{P}^{k-l-1}\setminus \mathfrak{P}^{l-k+1} \right) ~.
\] In addition $W_g g$ is continuous, hence $W_g g \in C_c(G)$. 

Finally, the implication $(b) \Rightarrow (a)$ is clear. 
\end{proof}

For later use we note that the proof of (c) $\Rightarrow$ (b) can be easily adapted to show the following statement. A converse to this observation will be shown in the following section.
\begin{corollary} \label{cor:wavelet_S0}
 Let $f,g \in \mathcal{S}_0(K)$. Then $W_g f \in C_c(G)$.

 More generally: Given $l \le k$ and $g \in \mathcal{S}_0(G)$, there exists a compact set $A \subset G$ such that for all $f \in \mathcal{S}_0(K) \cap S_l^k(K)$: ${\rm supp}(W_g f) \subset A$.
\end{corollary}

\begin{remark} \label{rem:explicit_wavelet}
 We stress that it is easy to construct functions $g \in \mathcal{S}_0(K)$ satisfying the admissibility condition from the lemma. As a particular case in point, we mention
 \[
  g = \mathbf{1}_{\mathfrak{D}} - q \mathbf{1}_{\mathfrak{P}}~.
 \]
\end{remark}

It turns out that the integration of wavelet coefficient functions $W_g f$ can be significantly simplified by the observation that these functions are locally constant whenever $g \in \mathcal{S}_0(K)$. Before we establish that fact, we introduce some notions related to compact open subgroups and their system of representatives.

\begin{definition} \label{defn:special}
 A subgroup $H<G = K \rtimes K^*$ is called \textbf{special compact open subgroup} if there exists an open compact subgroup $H_1 < (K,+)$ and a compact open subgroup $H_2 < (K^*,\cdot)$ such that
 \[
  H = \{ (x,h) :x \in H_1, h \in H_2 \}~.
 \]
Given a special compact open subgroup $H$, we call $R \subset G$ a \textbf{special system of representatives} modulo $H$ if
\[
 R = \{ (\lambda \gamma, \gamma): \lambda \in \Lambda, \gamma \in \Gamma \}
\]
where $\Gamma \subset K$ is a system of representatives modulo $H_1$, and $\Lambda \subset K^*$ a system of representatives modulo $H_2$.
\end{definition}

Special compact open subgroups $H$ and their special systems of representatives are particularly suitable for the discretization of weighted mixed $L^p$-norms.
The terminology regarding special systems of representatives is justified by the next lemma. The lemma also clarifies that special systems of representatives exist for every special compact open subgroup.
\begin{lemma} \label{lem:special_subgr}
Let $H_1 < (K,+)$ and $H_2 < (K^*,\cdot)$ be compact open subgroups.
\begin{enumerate}
 \item[(a)] The subset $H = \{ (x,h) :x \in H_1, h \in H_2 \} \subset G$ is a compact and open subset. It is a subgroup iff $H_1 \cdot H_2 \subset H_1$, where $\cdot$ denotes multiplication within the field $K$.
 \item[(b)] Assume that $H$ is a subgroup. If $\Gamma \subset K$ is a system of representatives modulo $H_1$, and $\Lambda \subset K^*$ a system of representatives modulo $H_2$, then
 \[
  R =  \{ (\lambda \gamma, \gamma): \lambda \in \Lambda, \gamma \in \Gamma \}
 \] is a system of representatives of the left $H$-cosets, i.e.,
 \[
  G = \bigcup_{r \in R}^\bullet r H~.
 \]
\end{enumerate}
\end{lemma}

\begin{proof}
 For the proof of (a), it is obvious from the semidirect product group law that $H_1 \cdot H_2 \subset H_1$ is sufficient to guarantee the subgroup property of $H$. Conversely, if $H$ is a subgroup, then
 \[
  (0,h) (x,1) = (h\cdot x,h) \in H
 \] whenever $h \in H_2$ and $x \in H_1$, and this entails $h\cdot x \in H_1$.

 For the proof of (b), let $(x,h) \in G$ be given. We are in search of unique tuple of elements $(\gamma,\lambda) \in \Gamma \times \Lambda$, $(x_1,h_1) \in H$ satisfying
 \[
  (x,h) = (\lambda \gamma, \gamma) (x_1, h_1) = (\gamma (\lambda + x_1) , \gamma h_1) ~.
 \]
 Comparing the second components on each sides, we can immediately identify a unique pair $(\gamma, h_1) \in \Gamma \times H_2$ with the property $h = \gamma h_1$.

 But then picking the unique pair $(\lambda, x_1)\in \Lambda \times H_1$ satisfying
 \[
  \lambda + x_1 = \gamma^{-1} x
 \]
 finishes the proof.
\end{proof}

We next identify a family of special open compact subgroups of particular relevance.
\begin{lemma} \label{lem:subgr_cosets}
 \begin{enumerate}
  \item[(a)] For $m \in \mathbb{N}$ let 
  \[
   \mathfrak{D}_m^* = \{ x \in K : |x-1|\leq q^{-m} \}~.
  \]
  Then $\mathfrak{D}_m^* \subset \mathfrak{D}^*$ is an open, compact multiplicative subgroup with index $[\mathfrak{D}^*: \mathfrak{D}_m^*] = q^m - q^{m-1}$.
  \item[(b)] For any $k \in \mathbb{Z}$ and $m \in \mathbb{N}$, the set $\mathfrak{P}^k \times \mathfrak{D}_m^* \subset G$ is a special open compact subgroup of $G$, with Haar measure $\mu_G(\mathfrak{P}^k \times \mathfrak{D}_m^*) = q^{-k-m}$.

  \item[(c)] For any $m \in \mathbb{N}$ let $\Lambda_m \subset \mathfrak{D}^*$ denote a system of representatives modulo $\mathfrak{D}_m^*$. Furthermore, let $\Gamma \subset K$ denote a system of representatives modulo the additive subgroup $\mathfrak{P}^k$. Then 
  \[
   R(m,\Lambda_m,\Gamma) = \{(\mathfrak{p}^n \lambda \gamma,  \mathfrak{p}^n \lambda): n \in \mathbb{Z},\lambda \in \Lambda_m, \gamma \in \Gamma \} \subset G
  \] is a special system of representatives modulo  $\mathfrak{P}^k \times \mathfrak{D}_m^* < G$.
 \end{enumerate}

\end{lemma}
\begin{proof}
 For the proof of $(a)$, first let $x \in \mathfrak{D}_m^*$, for $m \ge 1$. Since $|x-1| < |1|$, the ultrametric inequality gives
 \[
|x| \le \max(|1-x|,|1|) = 1 ~,
 \]
and thus $x \in \mathfrak{D}^*$.
 
We next show that $\mathfrak{D}_m^*$ is a multiplicative subgroup, for $m \ge 1$. For this purpose, let $1+x,1+y \in \mathfrak{D}_m^*$ be arbitrary elements, i.e., $x,y \in \mathfrak{P}^m$. Then the ultrametric triangle inequality entails
\[
 |(1+x)(1+y)-1| = |1+x+y+xy-1| \le \max(|x|,|y|,|xy|) \le p^{-m}~.
\] This shows $\mathfrak{D}_m^* \cdot \mathfrak{D}_m^* \subset \mathfrak{D}_m^*$. Furthermore, if $z = 1+x$, then the Neumann series
\[
 z^{-1} = \sum_{l=0}^\infty (-1)^l x^l
\] converges because of $|x^l| \leq q^{-ml}$, and because the field $K$ is complete. In fact, the choice of $x$ entails that the partial sums are contained in $\mathfrak{D}_m^*$, which is a closed subset of $K$. Hence the limit is contained in $\mathfrak{D}_m^*$, and $z^{-1} \in \mathfrak{D}_m^*$ follows.

For the computation of the index, we first observe that compactness of $\mathfrak{D}^*$ and openness of $\mathfrak{D}_m^*$ entail that the index is finite. Denoting a system of representatives mod $\mathfrak{D}_m^*$ in $\mathfrak{D}^*$ by $y_1,\ldots, y_r$, we have the disjoint union 
\[
 \mathfrak{D}^* = \bigcup_{j=1}^r y_j \mathfrak{D}_m^*~,
\] and since $|z| = 1$ holds on $\mathfrak{D}^* \supset y_j \mathfrak{D}_m^*$, we obtain for the multiplicative Haar measures
\[
 \int_{y_j \mathfrak{D}_m^*} \frac{dz}{|z|} = \int_{\mathfrak{D}_m^*} 
 dz = q^{-m}. 
\] This forces $r=q^m(1-q^{-1})= q^m - q^{m-1}$, as announced.

For the first part of $(b)$ we only need to check that $h x \in \mathfrak{P}^k$, for all $h \in \mathfrak{D}_m^*$ and all $x \in \mathfrak{P}^k$, which follows by
\[
 |hx| = |h| |x| = |x|. 
\] 
We then compute the Haar measure by the calculation
\begin{eqnarray*}
 \mu_G(\mathfrak{P}^k \times \mathfrak{D}_m^* \subset G) & = & 
 \int_{\mathfrak{P}^k} \int_{\mathfrak{D}_m^*} \frac{1}{|h|} dx dh \\
 & = & \int_{\mathfrak{P}^k} \int_{\mathfrak{D}_m^*}  dx dh \\
 & = & q^{-k-m}.
\end{eqnarray*}

In order to prove (c), it remains to show that
\[ \Lambda = \{ \mathfrak{p}^n \lambda :  n \in \mathbb{Z}, \lambda \in \Lambda_m \}
\] is a system of representatives modulo $\mathfrak{D}_m^*$ within $K^*$.
This follows from the fact that $\{ \mathfrak{p}^n : n \in \mathbb{Z} \}$ is a system of representatives modulo $\mathfrak{D}^*$ within $K^*$, and the assumption on $\Lambda_m$.
\end{proof}

Note that the family $(\mathfrak{P}^k \times \mathfrak{D}_m^*)_{m>0, k \in \mathbb{Z}}$ are a neighborhood basis at unity in $G$. As a consequence, every compact open subgroup $H < G$ contains a special compact open subgroup from this family, which is necessarily of finite index.

The full significance of the subgroups $\mathfrak{P}^k \times \mathfrak{D}_m^*$ is provided by the next observation, which shows that they arise naturally as fixed groups for the quasi-regular representation acting on functions $f \in \mathcal{S}^k_l(K)$. 

\begin{lemma}
 Let $g \in \mathcal{S}^k_l (K)$, and let $m = k-l \ge 0$. Then, for all $(x,h) \in \mathfrak{P}^k \times \mathfrak{D}_m^*$, one has 
 \[
  \pi(x,h) g = g~.
 \]
\end{lemma}
\begin{proof}
 By assumption, $g$ is supported on $\mathfrak{P}^l$ and constant on cosets of $\mathfrak{P}^k$. 
 
 Assume that $y \in K$ and $(x,h) \in \mathfrak{P}^k \times \mathfrak{D}_m^*$ are given. In the case where $h^{-1}(y-x) \in \mathfrak{P}^l$, or equivalently 
 \[
  |h^{-1}(y-x)|\leq q^{-l}~,
 \] the fact that $h \in \mathfrak{D_m}^*$ forces $y-x \in \mathfrak{P}^l$. Since we already assume $x \in \mathfrak{P}^k$, and the latter is a subgroup, this implies $y \in \mathfrak{P}^l$. 
 
 The contrapositive of this statement yields 
 \[
  y \not\in \mathfrak{P}^l \Rightarrow h^{-1}(y-x) \not\in \mathfrak{P}^l
 \]
 and it implies
 \[
 g(y) =  0 = g(h^{-1}(y-x)) = (\pi(x,h) g) (y) 
 \] whenever $y \in (K \setminus \mathfrak{P}^l)$. 

 On the other hand, whenever $y \in \mathfrak{P}^l$, the choice of $x \in \mathfrak{P}^k$ and $h \in \mathfrak{D}_m^*$ yields  
 \[
  |h^{-1}(y-x) - y| \le \max( |h^{-1}-1|~|y | , |x| ) \le \max(q^{-m-l},q^{-k}) = q^{-k}~.
 \] Since $g$ is constant on $\mathfrak{P}^k$-cosets, it again follows that 
 \[
  g(y) = g(h^{-1}(y-x)) = (\pi(x,h) g) (y) ~. 
 \]
\end{proof}

As a consequence, we note the following:
\begin{lemma} \label{lem:wgf_inv}
 Let $g \in \mathcal{S}^k_l(K)$, let $m=k-l$, and $f \in L^2(K)$. Then $W_g f$ is right-invariant under  $\mathfrak{P}^k \times \mathfrak{D}_m^*$, i.e., for all $s \in G$ and all $t \in \mathfrak{P}^k \times \mathfrak{D}_m^*$, 
 \[
  W_g f (st) = W_g f(s).
 \]
\end{lemma}
\begin{proof}
 This follows from the previous lemma via
 \[
  W_g f(st) = \langle f|\pi(st) g \rangle = \langle f| \pi(s) \pi(t) g \rangle = \langle f| \pi(s) g \rangle = W_g f(s)~. 
 \]
\end{proof}

The following lemma notes a complementary fact that will simplify the treatment of weighted $L^{p,q}$ norms. It is obtained by maximizing over cosets. Note that the resulting weight $w'$ is automatically continuous.
\begin{lemma}
Let $w$ be a weight on $G = K \rtimes K^*$, and let $H< K$ denote a compact open subgroup. Then there exists a weight $w' \sim w$ that is right-invariant under $H$.
\end{lemma}

\begin{lemma} \label{lem:integrate_cosets}
Let $H<G$ denote an open compact subgroup, and let $R \subset G$ denote a system of representatives for the left $H$-cosets.
Assume that $F: G \to \mathbb{C}$ is right $H$-invariant. Then
\begin{equation} \label{eqn:discr_integral}
 \int_G F(x) dx = \mu_G(H) \sum_{r \in R} F(r)~.
\end{equation} This implies in particular
\begin{equation} \label{eqn:discr_lp}
 \| F \|_{L^p} =  \mu_G(H)^{1/p} \| F|_{R} \|_{\ell^p}.
\end{equation}
Now assume that $0<s,t\leq\infty$ are given and let  $w: G \to \mathbb{R}^+$ denote a positive weight function on $G$ that is right $H$-invariant. Let $H = H_1 \times H_2 < G$ be a special open compact subgroup, and
\[
 R = \{ (\lambda \gamma, \lambda) ~:~ \lambda \in \Lambda, \gamma \in \Gamma \},
\] a special set of representatives, with suitable $\Gamma \subset K, \Lambda \subset K^*$.
Then
\begin{equation} \label{eqn:discr_mixed_norm}
 \| F \|_{L_w^{s,t}}^s = \mu_{K^*}(H_2) \mu_K(H_1)^{t/s} \sum_{\lambda \in \Lambda} |\lambda|^{t/s} \left( \sum_{\gamma \in \Gamma } |F(\lambda \gamma,\lambda)|^s w(\lambda \gamma,\lambda)^s \right)^{t/s}.
\end{equation}
The norm equalities hold in the extended sense that one side is finite iff the other side is.
\end{lemma}
\begin{proof}
 The proof of (\ref{eqn:discr_integral}) follows from the simple observation that
 \begin{eqnarray*}
  \int_G F(x) dx &  = &  \sum_{r \in R} \int_{r H} F(x) dx\\  & = & \sum_{r \in R} F(r) \mu_{G}(H)~.
 \end{eqnarray*}
 Applying this to $|F|^p$ yields (\ref{eqn:discr_lp}).

For the proof of (\ref{eqn:discr_mixed_norm}) first observe that the pointwise product of two right $H$-invariant functions (e.g. $F, w$) is again right $H$-invariant. Hence we may restrict the following computation to the case that $w$ is constant. For $s,t < \infty$ we have
\[
 \| F \|_{L^{s,t}}^t = \int_{K^*} \| F(\cdot, h) \|_s^t \frac{dh}{h} = \sum_{\lambda \in \Lambda} \int_{\lambda H_2} \| F(\cdot,h) \|_s^t \frac{dh}{h}~,
\]
where the second equality arises from the coset decomposition of $K^*$.
Given $h \in \lambda H_2$, we use the coset decomposition $K = \bigcup_{\gamma \in \Gamma} \lambda \gamma + \lambda H_1$ to compute
\[
\| F(\cdot, h) \|_s^s = \sum_{\gamma \in \Gamma} \int_{\lambda \gamma + \lambda H_1} |F(x,h)|^s dx = \sum_{\gamma \in \Gamma} \mu_K(\lambda H_1) |F(\lambda \gamma, \lambda)|^s~,
\] where the second equality used that if $x \in \lambda(\gamma + H_1)$ and $h \in \lambda H_2$, then $(x,h) \in (\lambda \gamma, \lambda) H$, and right $H$-invariance of $F$.

Note in particular the independence of this term of $h \in \lambda H_2$. Hence, summing over $\lambda$ and using $\mu_K(\lambda H_1) = |\lambda| \mu_K(H_1)$, one arrives at the desired equation.
\end{proof}

\section{Discrete wavelet expansions}

In this section, we aim to derive the existence of readily available tight wavelet frames arising from the action of a suitably chosen lattice acting on an arbitrary wavelet $0 \not= g \in \mathcal{S}_0(K)$. In principle, the methods of coorbit theory (which we only study in the subsequent section) could also be helpful here, however it turns out that for wavelets in $g \in \mathcal{S}_0(K)$, the arguments are drastically simplified by the fact that wavelet coefficients $W_g f$ are right-invariant under a compact open subgroup that only depends on $g$.

For the full appreciation of the following theorem, recall that Lemma \ref{lem:wgf_inv} applies to every element $g \in \mathcal{S}_0(G)$, and thereby allows to explicitly determine a compact open fixed group associated to a wavelet $g \in \mathcal{S}_0(G)$. In addition, the construction of associated systems of representatives is equally straightforward (see Lemma \ref{lem:subgr_cosets}). Hence the compact open subgroup $H$ and the system $R$ of representatives from the theorem are directly and explicitly computable from $g$. 

\begin{theorem} \label{thm:tight_frame}
 Let  $0 \not= g \in \mathcal{S}_0(K)$. Let $H < G$ be compact open with $\pi(x,h) g = g$ for all $(x,h) \in H$. Let $R \subset G$ be any system of representatives modulo $H$.
 Then the system $(\pi(r) g)_{ r \in R}$ is a tight frame of $L^2(K)$, satisfying for all $f \in L^2(K)$ the norm equality
   \begin{equation} \label{eqn:frame_l2}
    \sum_{r \in R} |\langle f| \pi(r) g \rangle|^2 = \frac{C_g}{\sqrt{\mu_G(H)}} \| f \|_2^2
   \end{equation}
   as well as the wavelet expansion
   \begin{equation} \label{eqn:framerec_l2}
    f = \frac{\sqrt{\mu_G(H)}}{C_g } \sum_{r \in R} \langle f| \pi(r) g \rangle \pi(r) g~.
   \end{equation}
\end{theorem}

\begin{proof}
 Equation (\ref{eqn:frame_l2}) follows from (\ref{eqn:adm_cond}) via Lemma \ref{lem:integrate_cosets}. The reconstruction formula (\ref{eqn:framerec_l2}) is a direct consequence.
\end{proof}

\begin{corollary} \label{cor:ortho_basis}
 Let $k\in \mathbb{N}$ and $H = \mathfrak{P}^k \times \mathfrak{D}_m^* < G, $ and $R \subset G$ be a special system of coset representatives modulo $H$ as in Lemma \ref{lem:subgr_cosets}. If $g = \check{\mathbf{1}}_{\mathfrak{p}^{-k}+\mathfrak{D}}$ then $\{\pi(r)g\,\,:\,\,r\in R\}$ forms a Parseval frame of $L^2(K)$. 
 \end{corollary}

 \begin{proof}
Since $\hat{g}$ is supported on $\mathfrak{P}^{-k}$ and constant on the cosets of $\mathfrak{D}$, $\hat{g} \in \mathcal{S}^0_{-k}$. So $g \in \mathcal{S}_0^k$ and $m = k-0 =k$. Also
$$C_g = \int_K \dfrac{|\hat{g(\xi)|^2}}{|\xi|}d\xi\,\,=\,\, \int_{\mathfrak{p}^{-k}+\mathfrak{D}} \dfrac{1}{|\xi|}d\xi\,\,=\,\, p^{-k}.$$
Since $\mu_G(\mathfrak{P}^k \times \mathfrak{D}_m^*) = p^{-k-m}=p^{-2k}$, we have   \begin{equation} \label{eqn:frame_l3}
    \sum_{r \in R} |\langle f| \pi(r) g \rangle|^2 =  \| f \|_2^2
   \end{equation}
   by Theorem \ref{thm:tight_frame}, which is the Parseval's identity.\par

\end{proof}

\begin{remark}
 Despite the extensive literature on the construction of wavelet orthonormal bases and frames on local fields, in particular in connection with multiresolution analysis, we are not aware of any source noting the rather simple approach towards the construction of wavelet frames outlined in Theorem \ref{thm:tight_frame}.
\end{remark}

As a consequence, we obtain wavelet characterizations and expansions for $\mathcal{S}_0(K)$ and its dual.

\begin{corollary} \label{cor:wavelet_char_S0}
Let $0 \not= g \in \mathcal{S}_0(K)$. Then, given any $f \in L^2(K)$, the following are equivalent:
\begin{enumerate}
 \item[(a)] $f \in \mathcal{S}_0(K)$.
 \item[(b)] $W_g f \in C_c(G)$.
 \item[(c)] Given any open compact subgroup $H<G$ and any set $R \subset G$ of coset representatives modulo $H$ with the property that $(\pi(r) g)_{r \in R}$ is a tight frame of $L^2(K)$ with frame constant $C$, the right-hand side of the frame reconstruction formula
 \[
  f = \frac{1}{C} \sum_{r \in R} W_g f(r) \pi(r) g
 \] is a finite sum.
\end{enumerate}
\end{corollary}
\begin{proof}
The implication (a) $\Rightarrow$ (b) was noted in Corollary \ref{cor:wavelet_S0}. For (b) $\Rightarrow$ (c) we note that by assumption on $g$, $H$ and $R$, the frame expansion holds for $f$ in the $L^2$-sense. In addition, assumption (b) yields that $W_g f$ is compactly supported, and the choice of $H$ and $R$ yields that the left translates of $H$ by $R$ yield an open, disjoint cover. This implies that only finitely many translates meet the support of $W_g f$, hence $W_g f(r) \not= 0$ holds only for finitely many $r \in R$.

For the implication (c) $\Rightarrow$ (a) observe that Theorem \ref{thm:tight_frame} yields the existence of $H$ and $R$ to which the assumption (c) can be applied. In particular, $f$ has a finite expansion in the wavelet system, the members of which are in $\mathcal{S}_0(K)$. Hence $f \in \mathcal{S}_0(K)$.
\end{proof}

For the characterization of the image of the discrete wavelet transform, it is useful to recall the fact that the image spaces of coefficient operators associated to tight frames are characterized by reproducing kernel relations. We want to apply the subsequent results to wavelet systems generated by the action of a suitable set $R$ on finitely many vectors $g_1,\ldots,g_m \in \mathcal{S}_0(K)$, such as wavelet ONBs constructed via multiresolution analyses. For this purpose, we introduce the notation $\underline{m} = \{ 1,\ldots, m \}$, for $m \in \mathbb{N}$. However one should bear in mind that the following results are all applicable (with $m=1$) to the tight wavelet frames provided by Theorem \ref{thm:tight_frame}.

Note that in all subsequent arguments, the open compact subgroup $H$ can be chosen independent of the mother wavelets $g_1,\ldots,g_m$; in particular, it is not required to be the fixed group of a wavelet $g_i$, say as in Theorem \ref{thm:tight_frame}. 

Note that from now on we shall use the notation $W^d$ for the \textit{discretized} wavelet transform associated to finitely many fixed wavelets $g_1,\ldots,g_m$.

\begin{lemma} \label{lem:kernel_S0}
 Let $g_1,\ldots,g_m \in \mathcal{S}_0(K)$, $H< G$ a compact open subgroup and $R \subset G$ a system of representatives mod $H$.

 Define
 \[
  (t_{(r,j),(s,k)})_{((r,j),(s,k) \in (R \times \underline{m}) \times (R \times \underline{m})} \in \mathbb{C}^{(R \times \underline{m}) \times (R \times \underline{m})}~,~ t_{(r,j),(s,k)} = \langle \pi(s) g_k | \pi(r) g_j \rangle~.
 \]
\begin{enumerate}
 \item[(a)] Letting
 \begin{equation} \label{eqn:rep_kernel}
  (T \alpha) (r,j) = \sum_{s \in R} \sum_{k=1}^m t_{(r,j),(s,k)} \alpha(s,k)
 \end{equation}  defines an operator $T: \mathbb{C}^{R \times \underline{m}} \to \mathbb{C}^{R \times \underline{m}}$. For any index $(r,j)$, the right-hand side of (\ref{eqn:rep_kernel}) is a finite sum. $T$ is continuous with respect to the product topology.
 \item[(b)] Let $\mathbb{C}^{(R \times \underline{m})} \subset \mathbb{C}^{\mathbb{R} \times \underline{m}}$ denote the subspace of finitely supported families. Then $T (\mathbb{C}^{(R \times \underline{m})}) \subset \mathbb{C}^{(R \times \underline{m})}$.
 \item[(c)] Assume in addition that $(\pi(r) g_j)_{r \in R, j \in \underline{m}}$ is a tight frame of $L^2(K)$ with frame constant $C>0$. Then $\frac{1}{C} T \circ \frac{1}{C} T = \frac{1}{C} T$. As a consequence, the subspace
 \[
  \mathbb{C}^{R \times \underline{m}}_T = \left\{ \alpha \in \mathbb{C}^{R \times \underline{m}} : \alpha = \frac{1}{C} T \alpha \right\}
 \]
 is closed in the product topology.
\end{enumerate}
\end{lemma}

\begin{proof}
First fix $(r,j) \in R \times \underline{m}$ and $k \in \underline{m}$. Then
\[
 \{ s\in R : \langle \pi(s) g_k  | \pi(r) g_j  \rangle \not= 0 \} \subset R \cap {\rm supp}(W_{g_k} \pi(r) g_j ) \}.
\]
As the $R$-shifts of $H$ provide a disjoint open cover of $G$, and the support of $W_{g_k} \pi(r) g_j$ is compact by Corollary \ref{cor:wavelet_S0}, the intersection is finite. But this entails that, for any fixed $(r,j) \in R \times \underline{m}$, the set
\[
 \{ (s,k) \in R \times \underline{m} : t_{(r,j),{(s,k)}}  \not= 0 \}
\] is finite. But this entails parts (a) and (b) of the Lemma.

For part (c), note that the tight frame assumption yields that the map $C^{-1/2} W^d: L^2(K) \to \ell^2(R \times \underline{m})$ is an isometry. This entails that $C^{-1} W^d \circ (W^d)^* : \ell^2(R \times \underline{m}) \to  \ell^2(R \times \underline{m}) $ is the orthogonal projection onto the closed subspace $W^d(L^2(K))$. Now plugging in definitions yields
\begin{eqnarray*}
 \left( C^{-1} W^d \circ (W^d)^* (\alpha) \right) (r,j)   & = & \left\langle C^{-1} \sum_{(s,k) \in R \times \underline{m}} \alpha(s,k) \pi(s) g_k| \pi(r) g_j \right\rangle \\
 & = & C^{-1} \sum_{(s,k) \in R \times \underline{m}} \alpha(s,k) \langle \pi(s) g_k| \pi(r) g_j \rangle \\
 & = & (T \alpha) (r,j)
 \end{eqnarray*}
showing that $P = C^{-1} T |_{\ell^2(R \times \underline{m})}$.

Note that $P^2 = P$ on $\ell^2(R \times \underline{m})$, and we want to use this to prove $C^{-2} T^2 = C^{-1} T$ on the larger space. To see this, note that by the proof of parts (a) and (b), for any index $(r,i)$ there exists a finite index set $M  \subset R \times \underline{m}$ such that the
$(T \alpha)(r,i)$ only depends on the restriction of $\alpha$ to $M$. By applying this argument to the elements of $M$, and taking the union, one finds a second finite set $M \subset M'  \subset R \times \underline{m}$ such that both $(T \alpha)(r,i)$ and $(T^2 \alpha)(r,i)$ only depend on the restriction of $\alpha$ to $M'$. Hence, letting $\beta \in \ell^2(R \times m)$ denote a finitely supported element coinciding on $M'$ with a given $\alpha \in \mathbb{C}^{R \times \underline{m}}$, we find that
\[
C^{-2} T^2 (\alpha) (r,i) = C^{-2} T^2 (\beta) (r,i) = P^2 \beta(r,i) = P (\beta)(r,i) = C^{-1} T(\alpha) (r,i)~.
\] Hence $C^{-1} T$ is an orthogonal projection, in particular its image space coincides with the kernel of the operator ${\rm Id} - C^{-1} T$. This operator is continuous, since $T$ is continuous, hence $\mathbb{C}_T^{R \times \underline{m}}$ is closed.
\end{proof}

We now discuss wavelet expansions of distributions $\varphi \in \mathcal{S}_0'(K)$. For this purpose, we extend the wavelet transform associated to $g \in \mathcal{S}_0(K)$ using the extended sesquilinear form
\[
 \forall \varphi \in \mathcal{S}_0'(K)~ \forall f \in \mathcal{S}_0(K)~:~\langle \varphi | f \rangle = \varphi(\overline{f})~.
\]
One then obtains a wavelet transform $\mathcal{S}_0'(K) \to G^\mathbb{C}$, that is continuous if one endows $\mathcal{S}_0'(K)$ with the weak-$*$-topology, and $G^\mathbb{C}$ with the topology of pointwise convergence.

\begin{proposition} \label{prop:expansion_S_strich}
 Let $g_1,\ldots,g_m \in \mathcal{S}_0(K)$, $H<G$ a compact open subgroup and $R \subset G$ a system of representatives mod $H$.
 \begin{enumerate}
  \item[(a)] Given any coefficient family $(\alpha(r,j))_{(r,j) \in R \times \underline{m}} \in \mathbb{C}^{R \times \underline{m}}$, the sum
  \[
   \sum_{(r,j) \in R \times \underline{m}} \alpha(r,j) \pi(r) g_j
  \] converges unconditionally in the weak-$*$-topology of $\mathcal{S}'_0(K)$. This defines a linear, continuous operator $S : \mathbb{C}^{R \times \underline{m}} \to \mathcal{S}_0'(K)$, when the left-hand side is endowed with the product topology, and the right-hand side carries the weak-$*$-topology.
  \item[(b)] Assume that $(\pi(r) g_j)_{r \in R}$ is a tight frame of $L^2(K)$, with frame constant $C$. Then the wavelet transform
 induces a bijection
 \[
  W^d : \mathcal{S}'_0(K) \to \mathbb{C}^{R \times \underline{m}}_T~, \varphi \mapsto (\langle  \varphi | \pi(r) g_j \rangle)_{(r,j) \in R \times \underline{m}}
 \]
 where $\mathbb{C}^{R \times \underline{m}}_T$ is defined with reference to the reproducing kernel operator $T$. If this space is endowed with the product topology, and $\mathcal{S}_0'(K)$ carries the weak-$*$-topology, then the discrete wavelet transform $W^d$ is a topological isomorphism.  The inverse of the operator is given by the wavelet inversion formula
 \[
  (\alpha(r,j))_{(r,j) \in R \times \underline{m}} \mapsto \frac{1}{C} \sum_{(r,j) \in R \times \underline{m}} \alpha(r,j) \pi(r) g_j~,
 \] thus entailing the unconditionally convergent wavelet expansion
 \begin{equation}
 \label{eqn:framerec_S0_dual}
  \varphi = \frac{1}{C} \sum_{(r,j) \in R \times \underline{m}} \langle \varphi | \pi(r) g_j \rangle \pi(r) g_j
 \end{equation}
 valid on $\mathcal{S}_0'(K)$.
 \end{enumerate}
\end{proposition}

\begin{proof}
For the proof of (a), we define a linear functional $\varphi$ by letting, for a given $f \in \mathcal{S}_0(K)$
\begin{equation}
  \label{eqn:defn_limit}
  \langle f| \varphi \rangle = \sum_{(r,j)} \langle f | \pi(r) g_j \rangle \overline{\alpha}_{r,j}~.
\end{equation}
For every $f \in \mathcal{S}_0(K)$, this sum has only finitely many non-zero terms (by Corollary \ref{cor:wavelet_S0}), hence it converges absolutely. But this means that $\varphi$ is a well-defined linear map, and the sum defining it converges unconditionally.

It is in fact continuous: For this purpose it is enough to recall that by Corollary \ref{cor:wavelet_S0}, for all $i=1,\ldots,m$, the supports of $W_{g_i} f$, for all $f \in \mathcal{S}_0(K) \cap S_l^k(K)$, are contained in a fixed compact set $A$, which is covered by finitely many translates of the open compact subgroup $H$. Accordingly, the right-hand side of (\ref{eqn:defn_limit}) has nonzero terms indexed by a fixed finite subset of the index set, for all $f \in  \mathcal{S}_0(K) \cap S_l^k(K)$, and the definition of the topology on $\mathcal{S}(K)$ then entails continuity with respect to the weak-*-topology. 




For part (b), we note that the map
\[
\frac{1}{C} S \circ W^d :  \varphi \mapsto \frac{1}{C} \sum_{(r,j) \in R \times \underline{m}} \langle \varphi| \pi(r) g_j \rangle \pi(r) g_j
\]
is well-defined, by part (a), and continuous with respect to the weak-$*$-topology.

Furthermore, it coincides with the identity map on $\mathcal{S}_0(K)$.
This subspace is weak-$*$-dense in $\mathcal{S}_0'(K)$. To see this, let $\varphi \in \mathcal{S}_0'(K)$, and let $U \subset \mathcal{S}_0'(K)$ denote a weak-$*$-neighborhood of $\varphi$. By definition of the topology, there exists $f_1,\ldots, f_k \in \mathcal{S}_0(K)$ and $\epsilon >0$ such that
\[
 \{ \psi \in \mathcal{S}_0'(K) : \forall j=1,\ldots, k ~|\langle f_k, \varphi \rangle - \langle f_k, \psi \rangle| < \epsilon \} \subset U.
\]
Now the Riesz-Fischer theorem applied to the space $\mathcal{H}_0 = \mathrm{span}\{f_1,\ldots,f_k \}$ provides an element $\varphi_0 \in \mathcal{H}_0$ such that
\[
 \forall j=1,\ldots, k~:~\langle f_j, \varphi_0 \rangle = \langle f_j, \varphi \rangle~.
\]
But that means $\varphi_0 \in \mathcal{S}_0(K) \cap U$.

Hence the density statement is shown, and the equation  $\frac{1}{C} S W^d \varphi = \varphi$ extends by continuity to all of $\mathcal{S}_0'(K)$, and we have proved (\ref{eqn:framerec_S0_dual}).

We next prove $ W^d \circ S = \frac{1}{C} T$.
Using unconditional convergence in the dual topology, we obtain
\begin{eqnarray*}
 W^d(S(\alpha))(r,j) & = & \frac{1}{C} \sum_{(s,k)} \langle \alpha(s,k) \pi(s) g_k | \pi(r) g_j \rangle
 \\ & = & \frac{1}{C} \sum_{(s,k)} \langle \pi(s) g_k | \pi(r) g_j \rangle \alpha(s,k)
 \\  & = & \frac{1}{C} T (\alpha) (r,j).
\end{eqnarray*}

In particular, since $\frac{1}{C} T \circ \frac{1}{C} T = \frac{1}{C} T$, we find that
$\frac{1}{C} W^d \circ S|_{\left(\mathcal{C}^{R \times \underline{m}}\right)_T} = {\rm Id}$.

In summary, $W^d: \mathcal{S}_0'(K) \to \left(\mathbb{C}^{R \times \underline{m}}\right)_T$  and $\frac{1}{C} S: \left(\mathbb{C}^{R \times \underline{m}}\right)_T \to  \mathcal{S}_0'(K)$ are continuous inverses of each other, which finishes the proof.
\end{proof}

\section{Wavelet coorbit theory on local fields}
 
 With the observations from the previous sections, the applicability of coorbit theory is now established with great ease.
 
 \begin{corollary} \label{cor:coorbit_appl}
  Let $g \in \mathcal{S}_0(K)$.
  Then $g \in \mathfrak{B}_v$, for any weight $v$.
 \end{corollary}
\begin{proof}
 By Lemma \ref{lem:g_adm} $W_g g \in C_c(G)$. For any relatively compact neighborhood $Q \subset G$ of the identity, this entails that $M_Q^R (W_gg)$ and $M_Q (W_g g)$ are also compactly supported and bounded. This guarantees weighted integrability of the maximal functions against any weight that is bounded on compact sets, and thus $g \in \mathcal{B}_v$. 
\end{proof}

\begin{remark}
 As a consequence, we find $\mathcal{S}_0(K) \subset \mathcal{H}^1_v$ for any continuous weight $v$. This inclusion is in fact continuous: For any sequence $f_n \to f$ converging in $\mathcal{S}_0(K)$, the fact that the supports of $f_n$ and are restricted to a fixed compact set produces by Lemma \ref{cor:wavelet_S0} that the wavelet transforms $W_g f_n$ are jointly supported in a fixed compact set, whenever $g \in \mathcal{S}_0(K)$. Since $f_n \to f$ in $L^2$ entails $W_g f_n \to W_g f$ uniformly, convergence in the coorbit space norm follows.

 Hence, the reservoir space $\mathcal{H}_v^{1,\sim}$ has a natural embedding in the antidual of $\mathcal{S}_0(K)$, which (after conjugation) coincides with the distribution space $\mathcal{S}_0'(K)$.

 Hence, it is convenient to consider wavelet coorbit spaces as subspaces of $\mathcal{S}_0
 '(K)$.
\end{remark}

The theory of Besov spaces over local fields has been developed significantly in recent years, extending many classical results from Euclidean settings to non-archimedean fields. Homogeneous Besov spaces over Vilenkin groups were introduced by W. Su and Onnweer in \cite{bib8}, where several characterisations were established including descriptions via mean oscillations and atomic decompositions. More recently in \cite{ashraf2023dilation}, these spaces were further explored focusing on the action of dilation operators and their localization properties.

In this section, we aim to provide a coorbit-theoretic characterization of certain homogeneous Besov spaces over local fields. Specifically, we demonstrate that these Besov spaces can be naturally realized as wavelet coorbit spaces over local fields, analogous to the well-established results over $\mathbb{R}$. We define homogeneous Besov spaces on a local field as following similarly to \cite{bib8}:
 \begin{definition}
     Let $\alpha \in \mathbb{R}$ and $0<s,t \leq \infty$, then
     $$\dot{B}_{\alpha, s,t} = \left\{ f \in \mathcal{S}_0'(K):\,\,\, \left\Vert f \right \Vert_{\dot{B}_{\alpha, s,t}} := \left(\sum_{k \in \mathbb{Z}} q^{k\alpha t}  \left\Vert f * \Phi_k \right \Vert_s^t \right )^{1/t} < \infty \right\}$$
      with usual modifications for $t = \infty$. Here we used $\Phi_k(x) = q^k \mathbf{1}_{\mathfrak{P}^k}(x) -  q^{k-1} \mathbf{1}_{\mathfrak{P}^{k-1}}(x)$.
 \end{definition}

 We can now identify homogeneous Besov spaces as coorbit spaces, as in the euclidean case (see \cite[§ 7.2]{bib2}).
 \begin{theorem} \label{thm:besov_coorbit}
  Let $1 \le s,t \le \infty$, and $\alpha \in \mathbb{R}$. Then  $\dot{B}_{\alpha, s,t} = Co(L^{s,t}_w)$, for the weight $w(x,h) = |h|^{-\alpha-1/2}$, with equivalent norms. 
 \end{theorem}
\begin{proof}
  Consider the wavelet $\phi(x) =  \mathbf{1}_{\mathfrak{D}}(x) -  q^{-1} \mathbf{1}_{\mathfrak{P}^{-1}}(x) $.  By construction $\phi \in \mathcal{S}_0(G)$. If $|h|=q^{-k}$ then
  \begin{eqnarray*}
      \phi(h^{-1}x) &=& \mathbf{1}_{\mathfrak{D}}(h^{-1}x) -   q^{-1}\mathbf{1}_{\mathfrak{P}^{-1}}(h^{-1}x)\\
      &=&   \mathbf{1}_{\mathfrak{P}^k}(x) -  q^{-1}\mathbf{1}_{\mathfrak{P}^{k-1}}(x) \\
      &=&  q^{-k} \Phi_{k}(x).
  \end{eqnarray*}
Note that $\Phi_k(-x) = \Phi_k(x)$ for all $x \in K$ as $\mathfrak{P}^k$ and $\mathfrak{P}^{k-1}$ are both subgroups of $K$. Now
\begin{eqnarray*}
  \langle f | \pi(x,h) \phi \rangle &=& |h|^{-1/2}\int_K f(t)\overline{\phi(h^{-1}(t-x))} dt\\
  &=& q^{k/2}\int_K f(t)(q^{-k}\Phi_{k}(t-x))dt\\
  &=& q^{-k/2}[ f * \Phi_{k}(x)].
\end{eqnarray*}
Now
 \begin{eqnarray*}
     \int_{K^*} \left \Vert w(\cdot,h) W_{\phi}f(\cdot,h)\right \Vert_s^t \dfrac{dh}{|h|} &=& \sum\limits_{k \in \mathbb{Z}} \int_{|h| = q^{-k}}  \left \Vert w(\cdot,h).W_{\phi}f(\cdot,h)\right \Vert_s^t \dfrac{dh}{|h|} \\
     &=& \sum\limits_{k \in \mathbb{Z}} \int_{|h| = q^{-k}} |w(h)|^t q^{(-\frac{kt}{2})}  \left\Vert f * \Phi_{k} \right \Vert_s^t \dfrac{dh}{|h|}.
 \end{eqnarray*}

 Hence we see that taking $w(h) = |h|^{-\alpha - 1/2}$ allows to actually conclude the \textit{equality} $$ \left \Vert f \right \Vert_{Co(L_w^{s,t})} = \left\Vert f \right \Vert_{\dot{B}_{\alpha, s,t}},$$ if the wavelet $\phi$ is used for the definition of the coorbit space norm.
\end{proof}

We now turn to frame characterizations of the coorbit spaces $Co(Y)$. This requires the introduction of suitable norms in coefficient spaces.

\begin{defrem} \label{defrem:coeff_space}
 Let $R \subset G$ denote a system of representatives of the left cosets modulo an open compact subgroup $H$. Then $R \subset G$ is a well-spread set in the terminology of \cite{MR1021139}.

 Now let $Y$ denote a solid BF-space on $G$. Following \cite[Definition 3.4]{MR1021139}, we define the discrete coefficient space norm associated to $Y$ as
 \begin{equation} \label{eqn:norm_Yd}
  \| \alpha \|_{Y_d} = \left\| \sum_{r \in R} \alpha_r \mathbf{1}_{rH} \right\|_{Y}
 \end{equation}
and the associated coefficient space as
\begin{equation} \label{eqn:Yd}
Y_d = \left\{ \alpha \in \mathbb{C}^R : \| \alpha \|_{Y_d} < \infty \right\}~.
\end{equation}

Since we want to cover multi-wavelet systems of the type $(\pi(r) g_j)_{(r,j) \in R \times \underline{m}}$, we define for $\alpha \in \mathbb{C}^{R \times \underline{m}}$
\[
 \left\| \alpha \right\|_{Y_d} = \left( \sum_{j=1}^m \| \alpha(\cdot, j) \|_{Y_d}^2 \right)^{1/2}~.
\] and define the space $Y_d \subset \mathbb{C}^{R \times \underline{m}}$ accordingly.
\end{defrem}

\begin{remark}
 For $Y = L^{s,t}_w(G)$ and $m=1$, it is often possible to give more explicit meaning to the norm $\| \cdot \|_{Y_d}$.  E.g., in the cases $s=t$, or when $R$ and $H$ are special in the sense of Definition \ref{defn:special}, Lemma \ref{lem:integrate_cosets} implies
\[
 \| \alpha \|_{Y_d} = \| \alpha \|_{\ell^{s,t}_{w'}(R)}~.
\]
for the weight $w'(x,h) = |h|^{t/s} w(x,h)$.
\end{remark}

The following theorem can then be seen as a summary of the frame properties of discrete wavelet systems in various spaces.
\begin{theorem} \label{thm:disc_char_coorbit}
 Let $g_1,\ldots,g_m \in \mathcal{S}_0(K)$, $H<G$ a compact open subgroup and $R \subset G$ a system of representatives mod $H$. Assume that the system $(\pi(r) g_j)_{(r,j) \in R \times \underline{m}}$ is a tight wavelet frame of $L^2(K)$.

 Fix $X \in \{ \mathcal{S}_0 (K), \mathcal{S}_0'(K),Co(Y) \}$, for a solid BF space $Y$. Denote by $X_d \subset \mathbb{C}^{R \times \underline{m}}$ the associated coefficient space, i.e.,
 \[
  X_d = \mathbb{C}^{(R \times \underline{m}) } \mbox{ if } Y = \mathcal{S}_0(K)~,~ X_d =  \mathbb{C}^{R \times \underline{m} } \mbox{ if } Y =
  \mathcal{S}_0'(K)~ \] and \[ X_d = Y_d \mbox{ if } X = Co(Y)~.
 \]
 Furthermore, we let $(X_d)_T = X_d \cap \mathbb{C}^{R \times \underline{m}}_T$, the elements of $X_d$ contained in the reproducing kernel space associated to the tight wavelet frame via Lemma \ref{lem:kernel_S0}.

 Given $f \in \mathcal{S}_0'(K)$, let $W^d(f) = \left( \langle f| \pi(r) g_j \rangle \right)_{(r,j) \in R \times \underline{m}}$.
 Then the following statements hold:
 \begin{enumerate}
  \item[(a)] For every $f \in \mathcal{S}_0'(K)$: $f \in X$ iff $W^d(f) \in X_d$.
  \item[(b)] \textbf{Banach frame property:} The discrete wavelet transform $W^d: X \to (X_d)_T$ is a bijection. For all $Y$ except $\mathcal{S}_0(K)$, $W^d$ is topological. In particular, for the coorbit space case $\| f \|_X \asymp \| W^d(f) \|_{X_d}$.
  \item[(c)] \textbf{Atomic decomposition:} Assume that the finite sequences are dense in $X_d$, or that $X = \mathcal{S}_0'(K)$. Then every $f \in X$ has the wavelet expansion
  \[
   f = \sum_{(r,j) \in R \times \mathbb{R}} \langle f | \pi(r) g_j \rangle \pi(r) g_j~,
  \] converging unconditionally in the topology on $X$.
  \item[(d)] Assume that the wavelet system is an orthonormal basis. Then $(X_d)_T = X_d$. As a consequence the discrete wavelet transform is a bijection onto all of $X_d$. The wavelet system is simultaneously an unconditional basis of all coorbit spaces $Y$ with the property that the finitely supported sequences are dense in $Y_d$. It is furthermore an algebraic basis of $\mathcal{S}_0(K)$: Every $f \in \mathcal{S}_0(K)$ is a unique linear combination of elements from the wavelet basis.
 \end{enumerate}

\end{theorem}

\begin{proof}
 For the cases $Y \in \{ \mathcal{S}_0(K), \mathcal{S}_0'(K) \}$, the statements (a)-(c) are already settled, by Lemma \ref{lem:kernel_S0} and Proposition \ref{prop:expansion_S_strich}. Furthermore, the wavelet system is an orthonormal basis iff the frame constant equals $C=1$, and the kernel $T$ is the identity matrix. This entails $\mathbb{C}^{R \times \underline{m}}_T = \mathbb{C}^{R \times \underline{m}}$, and thus part (d) for $Y \in \{ \mathcal{S}_0(K), \mathcal{S}_0'(K) \}$.

 It remains to consider the coorbit space case.
 We first show the continuity of the wavelet transform $W^d: X \to X_d$. By definition of the norm on $X_d$, we can restrict to the case $m=1$, which is settled by Theorem 6.1(i) of \cite{MR1021139}.

 Now consider the synthesis operator $S: X_d \to \mathcal{S}_0'(K)$
 \begin{equation} \label{eqn:synth_op}
  S(\alpha) = \sum_{(r,j) \in R \times \underline{m}} \alpha(r,j) \pi(r) g_j~.
 \end{equation}
By Proposition \ref{prop:expansion_S_strich}, the right-hand side converges unconditionally in $\mathcal{S}_0'(K)$.
We claim that $S$ is a bounded operator into $X = Co(Y)$. Clearly it is enough to prove this for the case of $m=1$, and we let $g = g_1$ accordingly. Convergence in $\mathcal{S}_0'(K)$ translates to pointwise convergence of wavelet coefficients. Hence, using the covariance properties of the wavelet transform we obtain for $\varphi = S(\alpha) \in \mathcal{S}_0'(K)$ that
\[
 W_g \varphi = \sum_{r \in R} \alpha(r) \lambda(r) W_g g~.
\] where $\lambda(r)$ denotes the left shift by $r$. Now the right-hand side converges unconditionally in $Y$ by Lemma 2.4.13 of \cite{bib5}, and this entails unconditional convergence of the right-hand side in the coorbit space norm: Given $\epsilon >0$, there exists a finite set $S \subset \mathbb{R}$ such that for all finite $S'$ with $S \subset S' \subset R$:
\[
\epsilon > \left\| \sum_{r \in S'}  \alpha(r) \lambda(r) W_g g - \sum_{r \in S}  \alpha(r) \lambda(r) W_g g \right\|_Y = \left\| W_g \left(\sum_{r \in S' \setminus S} \alpha(r) \pi(r) g \right) \right\|_Y ~.
\]
But by definition of coorbit space norms, this amounts to
\[
\left\| \sum_{ r \in S' \setminus S} \alpha(r) \pi(r) g \right\|_{Co Y} < \epsilon~,
\] and unconditional convergence in $\| \cdot \|_X$ is established.

By Corollary \ref{cor:coorbit_appl}, $W_g g \in W^R(L^\infty,L^1_v)$ for a control weight $v$ associated to $Y$. Hence Proposition 5.2 of \cite{MR1021139} (applicable whenever the finite sequences are dense in $X_d$) provides the estimate
\[
\| \varphi \|_{Co(Y)} = \| W_g \varphi \|_{Y} \le C \| (\alpha_r) \|_{Y_d}~,
\] which shows that $S$ is bounded.

Furthermore, we know $W^d : X \to (X_d)_T$ by Proposition \ref{prop:expansion_S_strich}, and $W^d$ and $\frac{1}{C} S$ are two-sided inverses of each other. Hence $W^d: X \to X_d$ is a topological isomorphism.

 Hence (a)-(c) are shown for the coorbit case. Part (d) follows again from $(X_d)_T = X_d$, if the wavelet system is an orthonormal basis.
\end{proof}

The next two remarks exhibit explicit examples of discrete wavelet systems to which Theorem \ref{thm:disc_char_coorbit} applies. We begin with the tight frame case: 
\begin{remark} \label{rem:simple_frame} A simple wavelet is provided by Remark \ref{rem:explicit_wavelet}, namely 
\[
g = \mathbf{1}_{\mathfrak{D}} - q \mathbf{1}_{\mathfrak{P}}~.
\] $g$ is supported in $\mathfrak{D} $, and constant on cosets mod $\mathfrak{P}$, i.e., $g \in \mathcal{S}_{0}^{1}$. This entails that $W_g f$ is invariant under $\mathfrak{P} \times \mathfrak{D}_1^*$, by Lemma \ref{lem:wgf_inv}. Hence picking a system of representatives
\[
\Lambda = \{ \lambda_1,\ldots,\lambda_{q-1} \} \subset \mathfrak{D}^*
\] of $\mathfrak{D}_1^* \subset \mathfrak{D}^*$, as well as a system $\Gamma \subset K^+$ of representatives modulo $\mathfrak{P}$ allows to conclude via Theorem \ref{thm:tight_frame} that the wavelet system  
\[
\left( \pi(r) g \right)_{r \in R}~,~
\] 
based on 
\[
 R = \{ \left( \mathfrak{p}^n \lambda_j \gamma, \mathfrak{p}^n \lambda_j \right) : k \in \mathbb{Z}, j = 1,\ldots, q-1, \gamma \in \Gamma \}~
\]
is a tight frame, that has the further properties formulated in Theorem \ref{thm:disc_char_coorbit} (a)-(c). In particular, the system is a Banach frame, giving rise to explicit atomic decompositions converging simultaneously in all relevant coorbit space norms.  
\end{remark}

\begin{remark} \label{rem:ex_wavelet_onb}
We expect that many of the existing constructions of wavelet orthonormal bases on local fields fall in the realm of Theorem \ref{thm:disc_char_coorbit}, at the expense of checking additional fineprint. 

As an example, we give an explicit construction of a wavelet orthonormal basis with several generators from the literature, and show how Theorem \ref{thm:disc_char_coorbit} applies. To a large extent, the following relies on the constructions of wavelet ONBs in \cite{bib9,bib10}, in particular on Theorems 4.2 and 4.3 of \cite{bib10}. It will turn out  that the translation operators used to define wavelet systems in \cite{bib9,bib10} are somewhat different from the ones used here. Ultimately, this turns out to be a minor obstacle, but it is probably indicative of the problems one may encounter in the study of other discrete wavelet systems on local fields. 

We now present the construction of a wavelet system following \cite[Theorem 4.2]{bib10}, by systematically defining the various objects in said result, using (for the most part) the same notation as in \cite{bib10}. 

To begin with, we fix the open compact subgroup $H = \mathfrak{D}< K^+$. Then, using a character $\chi$ with ${\rm Ker}(\chi) = \mathfrak{D}$ to identify $K^+$ with its dual group, we find that the annihilator subgroup of $H^\bot$ is $H$ itself. The map $A: x \mapsto \mathfrak{p}^{-1}x$ defines an automorphism of $K^+$, and it is straightforward to check that 
$A$ is indeed expansive with respect to $H$, in the sense of \cite[Definition 2.2]{bib10}. Furthermore, under our identification of $K^+$ with its dual via the fixed character $\chi$, the dual automorphism $A^*$ coincides with $A$, and $A^* H^\bot = \mathfrak{p}^{-1} \mathfrak{D}$. 

We next pick representatives $\sigma_0,\ldots, \sigma_{q-1}$ of $\mathfrak{p}^{-1} \mathfrak{D}/\mathfrak{D}$, with the convention that $\sigma_0 = 0 \in K$, and define the mother wavelets $g_i$, for $i=1,\ldots,q-1$ via the Fourier transform, namely
\[
\widehat{g}_i = \mathbf{1}_{\sigma_i + \mathfrak{D}}.
\]
Note that since $0 \not\in {\rm supp}(\widehat{g}_i)$, we have $g_i \in \mathcal{S}_0(K)$.

Now Theorem 4.2 of \cite{bib10} implies that the system 
\begin{equation} 
\label{eqn:onb_benedetto}
\{ \delta_A^n \tau_{[s]} g_i : 1 \le i \le q-1, n \in \mathbb{Z}, [s] \in K^+/H \}
\end{equation} is an orthonormal basis of $L^2(K)$, and it remains to establish how this system relates to the wavelet constructions treated in Theorem \ref{thm:disc_char_coorbit}.

Firstly, it is immediate that the dilation operator $\delta_A$ coincides with the unitary dilations implemented by the quasiregular representation, i.e., 
\[
\delta_A = \pi(0,\mathfrak{p}^{-1})~. 
\] 
Thus it remains to compare the translation operators $\tau_{[s]}$ with $\pi(s,0)$. These are in fact different; in particular, $\tau_{[s]}$ is $H$-invariant by definition, i.e., $\tau_{[sh]} = \tau_{[s]}$ whenever $h \in H$, whereas this property clearly does not hold for the translation operators. (Note that this property of the $\tau_{[s]}$ is the chief justification for using the quotient group $K^+/H$ as index set for the translations in \cite{bib10}.)

Now Theorem 4.3 of \cite{bib10} provides a formula for $\tau_{[s]} g_i$, which our particular identification of $K^+$ with $\hat{K}^+$ allows to write as 
\[
(\tau_{[s]} g_i) = \chi(\sigma_i x) \mathbf{1}_{s+H}(x)~,
\] for any $s \in G$. Comparing this expression for $s = 0$ and arbitrary $s \in K$ then gives 
\[
(\pi(s,1) g_i)(x) = g_i(x-s) = \chi(\sigma_i (x)-s) \mathbf{1}_{H}(x-s) = \overline{\chi(\sigma_i s)} \left(\tau_{[s]} g_i \right)(x)~,
\]
i.e., 
\[
\pi (s,1) g_i = \overline{\chi(\sigma_i s)} \tau_{[s]} g_i~,
\]
and finally
\[
\pi (\mathfrak{p}^n s, \mathfrak{p}^n) g_i = \pi( 0, \mathfrak{p}^n) \pi(s,1) g_i =  \overline{\chi(\sigma_i s)} ~ \delta_A^{-n} \tau_{[s]} g_i.
\]
But this entails that if we now fix a system of representatives $\Gamma \subset K$ modulo $H$, and define the subset
\[
R = \left\{ (\mathfrak{p}^n \gamma, \mathfrak{p}^n) : n \in \mathbb{Z}, \gamma \in \Gamma \right\}
\] then $R$ is a system of representative in $K \rtimes K^*$ modulo the open compact subgroup $\mathfrak{D} \times \mathfrak{D}^*$, and the system 
\[
\left( \pi(r) g_i \right)_{i =1,\ldots,q-1, r \in R}
\] coincides with the wavelet system (\ref{eqn:onb_benedetto}), except for the choice of index set and unimodular factors. Hence it is an orthonormal basis as well, and part (d) of Theorem \ref{thm:disc_char_coorbit} applies to it, establishing that this system is a joint unconditional basis of the full scale of coorbit spaces. 

We note in passing that the coefficient spaces associated to solid function spaces $Y$ via Definition \ref{defrem:coeff_space} are again solid, and in particular invariant under pointwise multiplication with factors of modulus one. Hence the wavelet basis (\ref{eqn:onb_benedetto}), with its particular choice of normalizations, is indeed also an unconditional basis of the scale of coorbit spaces. 
\end{remark}

The following remark collects the different facets of wavelet analysis that drastically simplify over local fields, whenever wavelets $g \in \mathcal{S}^k_l$ are employed.  Arguably, the main structural reason for this phenomenon is the fact that all groups $K,K^*$ and $G = K \rtimes K^*$ have a neighborhood basis of unity consisting of compact open subgroups.

\begin{remark}
 \begin{enumerate}
  \item The observation that $W_g g \in C_c(G)$ holds for all wavelets in $g \in \mathcal{S}_k^l(K)$ is in stark contrast to the euclidean setting, where a well-known \textbf{qualitative uncertainty principle} due to Wilczok \cite{MR1758876} states that $W_g f$ is compactly supported only if either $f= 0$ or $g=0$.
  \item The view of wavelet frames as discretizations of continuous wavelet systems provides a reasonable motivation and orientation for the construction of such systems. Arguably the development of coorbit space theory can be seen as a rather faithful implementation of this idea, and the techniques developed to make this approach work (oscillation estimates and Wiener amalgam spaces) illustrate the considerable difficulties one encounters in the implementation over generally locally compact groups.

  In particular, frame construction following this approach usually requires a detailed study of density conditions and related frame bounds, and does not yield tight frames. As a consequence, dual frames exist in principle, but may not have the structure of wavelet frames themselves; they are therefore generally hard to determine explicitly. 
  
  Furthermore, explicitly determining candidate wavelets and associated sampling densities, for which the approach yields the desired sampling results is anything but trivial. For wavelet systems in euclidean domains one typically employs vanishing moment conditions for such a purpose; see e.g. \cite{MR2327472,MR2898467,MR3452925}. As these sources illustrate, this leads to rather technical, and typically very pessimistic assumptions. 
  
  Tight wavelet frame constructions in the euclidean setting rest on the development of dedicated, Fourier-analytic criteria, using notions such as multiresolution analysis. See e.g. the classical sources \cite{MR0836025,MR1162107} which inspired a large body of subsequent work in this direction, including adaptations to local fields.
  
  It is notable that all of these complications vanish in the local field case. By comparison to the methods for the euclidean case, but also by comparison to the construction of muliresolution analyses, the construction of tight wavelet frames following Theorem \ref{thm:tight_frame} is almost effortless: Any nonzero $g \in \mathcal{S}_0(K)$ can serve as a wavelet. Given such a wavelet, the sampling density is easily determined (Lemma \ref{lem:wgf_inv}), as is the sampling set, and one automatically obtains a tight frame with associated explicit atomic decompositions valid for all coorbit spaces.  
   \item The very clean and lossless sampling theory for wavelet transforms is enabled by the fact that the sampling set $R$ in question can be chosen a set of coset representatives modulo an open, compact subgroup, which acts as a fixed group under the quasi-regular representation. This aspect is another feature of the local field setting that does not have an analog in the Euclidean setting: The affine group of the reals only has noncompact closed subgroups. As a consequence, one can show that the quasiregular action of the affine group on any nonzero square-integrable function is necessarily free. 
  \item The analogous remarks hold for the discretization of more general coorbit space norms. Again, this is all greatly facilitated by the properties of functions in $\mathcal{S}_0(K)$, leading to stronger results and easier proofs.
\item While our presentation of coorbit space theory over the affine group of the local field $K$ has been largely self-contained up to Theorem \ref{thm:disc_char_coorbit}, in the proof of that result we did resort to citing key convergence statements from the literature. In principle, we expect that an adaptation of the cited results to our setting is possible, and that the adapted proofs may again be substantially simpler than the proofs for the fully general case considered in the original sources. For example, the norm equivalence statement from Theorem \ref{thm:disc_char_coorbit} (b) is more or less directly available via Lemma \ref{lem:integrate_cosets} , at least for coorbit spaces associated to the spaces $L^{p,q}_v(G)$ (which include the Besov spaces). 

However, the chief purpose of this paper is to provide easy access to coorbit theory for the special setting of wavelets over local fields, rather than to provide a fully self-contained account.
 \end{enumerate}
\end{remark}

We finally remark on possible extensions of the results in this paper:
\begin{remark}
 \begin{enumerate}
  \item This paper concentrated on the Banach space setting, and avoided the quasi-Banach setting. We expect that the arguments largely go through also for coorbit spaces associated to quasi-Banach spaces, and a strong motivation to establish this extension is provided by the homogeneous Besov space $\dot{B}_{s,t,\alpha}(K)$ with either $0< s < 1$ or $0<t<1$; these are not covered by our present results. The relevant source for this problem will be the paper \cite{van2022coorbit}.
  \item One might also consider higher-dimensional wavelet systems on $L^2(K^n \rtimes H)$, with a suitably matrix group $H < Gl(n,K)$. The coorbit theory for the euclidean setting is fairly well-developed (see e.g. \cite{bib6,bib7}), and it might be interesting to see whether the great simplifications for local fields carry over to the higher-dimensional systems.
 \end{enumerate}

\end{remark}



\bibliographystyle{abbrv}
\bibliography{sn-bibliography-2}

\begin{thebibliography}{10}

\bibitem{MR2673705}
S.~Albeverio, S.~Evdokimov, and M.~Skopina.
\newblock {$p$}-adic multiresolution analysis and wavelet frames.
\newblock {\em J. Fourier Anal. Appl.}, 16(5):693--714, 2010.

\bibitem{ashraf2023dilation}
S.~Ashraf and Q.~Jahan.
\newblock Dilation operators in besov spaces over local fields.
\newblock {\em Advances in Operator Theory}, 8(2):27, 2023.

\bibitem{mabook}
B.~Behera and Q.~Jahan.
\newblock {\em Wavelet Analysis on Local Fields of Positive Characteristic}.
\newblock Springer Singapore, 2021.

\bibitem{bib9}
J.~J. Benedetto and R.~L. Benedetto.
\newblock A wavelet theory for local fields and related groups.
\newblock {\em The Journal of Geometric Analysis}, 14:423–456, 2004.

\bibitem{bib10}
R.~L. Benedetto.
\newblock Examples of wavelets for local fields.
\newblock In {\em Wavelets, frames and operator theory}, volume 345 of {\em Contemp. Math.}, pages 27--47. Amer. Math. Soc., Providence, RI, 2004.

\bibitem{MR1162107}
I.~Daubechies.
\newblock {\em Ten lectures on wavelets}, volume~61 of {\em CBMS-NSF Regional Conference Series in Applied Mathematics}.
\newblock Society for Industrial and Applied Mathematics (SIAM), Philadelphia, PA, 1992.

\bibitem{MR0836025}
I.~Daubechies, A.~Grossmann, and Y.~Meyer.
\newblock Painless nonorthogonal expansions.
\newblock {\em J. Math. Phys.}, 27(5):1271--1283, 1986.

\bibitem{Far}
Y.~A. Farkov.
\newblock Multiresolution analysis and wavelets on vilenkin groups.
\newblock {\em Facta Universitatis (NIS), series: Electrical Energy}, 21:309--325, 2008.

\bibitem{bib2}
H.~G. Feichtinger and K.~Gr{\"o}chenig.
\newblock A unified approach to atomic decompositions via integrable group representations.
\newblock In {\em Function Spaces and Applications}, pages 52--73, Berlin, Heidelberg, 1988. Springer Berlin Heidelberg.

\bibitem{MR1021139}
H.~G. Feichtinger and K.~Gr\"ochenig.
\newblock Banach spaces related to integrable group representations and their atomic decompositions. {I}.
\newblock {\em J. Funct. Anal.}, 86(2):307--340, 1989.

\bibitem{MR2327472}
H.~G. Feichtinger, W.~Sun, and X.~Zhou.
\newblock Two {B}anach spaces of atoms for stable wavelet frame expansions.
\newblock {\em J. Approx. Theory}, 146(1):28--70, 2007.

\bibitem{MR2130226}
H.~F\"{u}hr.
\newblock {\em Abstract harmonic analysis of continuous wavelet transforms}, volume 1863 of {\em Lecture Notes in Mathematics}.
\newblock Springer-Verlag, Berlin, 2005.

\bibitem{bib4}
H.~F{\"u}hr.
\newblock Generalized {C}alder{\'o}n conditions and regular orbit spaces.
\newblock In {\em Colloquium Mathematicum}, volume 120, pages 103--126. Institute of Mathematics Polish Academy of Sciences, 2010.

\bibitem{bib6}
H.~F{\"u}hr.
\newblock Coorbit spaces and wavelet coefficient decay over general dilation groups.
\newblock {\em Transactions of the American Mathematical Society}, 367(10):7373--7401, 2015.

\bibitem{MR3452925}
H.~F\"uhr.
\newblock Vanishing moment conditions for wavelet atoms in higher dimensions.
\newblock {\em Adv. Comput. Math.}, 42(1):127--153, 2016.

\bibitem{bib7}
H.~F{\"u}hr and F.~Voigtlaender.
\newblock Wavelet coorbit spaces viewed as decomposition spaces.
\newblock {\em Journal of Functional Analysis}, 269(1):80--154, 2015.

\bibitem{bib3}
K.~Gr{\"o}chenig.
\newblock Describing functions: atomic decompositions versus frames.
\newblock {\em Monatshefte f{\"u}r Mathematik}, 112:1--42, 1991.

\bibitem{MR2558153}
A.~Y. Khrennikov, V.~M. Shelkovich, and M.~Skopina.
\newblock {$p$}-adic refinable functions and {MRA}-based wavelets.
\newblock {\em J. Approx. Theory}, 161(1):226--238, 2009.

\bibitem{MR1918846}
S.~V. Kozyrev.
\newblock Wavelet theory as {$p$}-adic spectral analysis.
\newblock {\em Izv. Ross. Akad. Nauk Ser. Mat.}, 66(2):149--158, 2002.

\bibitem{MR1373159}
W.~C. Lang.
\newblock Orthogonal wavelets on the {C}antor dyadic group.
\newblock {\em SIAM J. Math. Anal.}, 27(1):305--312, 1996.

\bibitem{bib8}
C.~Onneweer and S.~Weiyi.
\newblock Homogeneous {B}esov spaces on locally compact {V}ilenkin groups.
\newblock {\em Studia Mathematica}, 93(1):17--39, 1989.

\bibitem{MR2511868}
V.~Shelkovich and M.~Skopina.
\newblock {$p$}-adic {H}aar multiresolution analysis and pseudo-differential operators.
\newblock {\em J. Fourier Anal. Appl.}, 15(3):366--393, 2009.

\bibitem{taible}
M.~H. Taibleson.
\newblock {\em Fourier Analysis on Local Fields}, volume~15.
\newblock Princeton University Press, 2015.

\bibitem{MR2898467}
T.~Ullrich.
\newblock Continuous characterizations of {B}esov-{L}izorkin-{T}riebel spaces and new interpretations as coorbits.
\newblock {\em J. Funct. Spaces Appl.}, pages Art. ID 163213, 47, 2012.

\bibitem{van2022coorbit}
J.~T. van Velthoven and F.~Voigtlaender.
\newblock Coorbit spaces associated to quasi-banach function spaces and their molecular decomposition.
\newblock {\em arXiv preprint arXiv:2203.07959}, 2022.

\bibitem{bib5}
F.~Voigtlaender.
\newblock {\em Embedding theorems for decomposition spaces with applications to wavelet coorbit spaces}.
\newblock PhD thesis, RWTH Aachen, 2016.

\bibitem{MR1758876}
E.~Wilczok.
\newblock New uncertainty principles for the continuous {G}abor transform and the continuous wavelet transform.
\newblock {\em Doc. Math.}, 5:201--226, 2000.

\end{thebibliography}



\end{document}